\documentclass{amsart}
\usepackage{amssymb}
\usepackage{amscd}

\newtheorem{theorem}{Theorem}
\newtheorem{proposition}[theorem]{Proposition}
\newtheorem{definition}[theorem]{Definition}

\newtheorem{question}[theorem]{Question}
\newtheorem{lemma}[theorem]{Lemma}

\numberwithin{equation}{section}
\numberwithin{theorem}{section}

\begin{document}
\title[$L^p$ isometries]{On the structure of isometries between noncommutative $L^p$ spaces}
\author{David Sherman}
\address{Department of Mathematics\\ University of Illinois\\ Urbana, IL 61801}
\email{dasherma@math.uiuc.edu}

\begin{abstract}
We prove some structure results for isometries between noncommutative $L^p$ spaces associated to von Neumann algebras.  We find that an isometry $T: L^p(\mathcal{M}_1) \to L^p(\mathcal{M}_2)$ ($1 \le p < \infty$, $p \ne 2$) can be canonically expressed in a certain simple form whenever $\mathcal{M}_1$ has variants of Watanabe's extension property [W2].  Conversely, this form always defines an isometry provided that $\mathcal{M}_1$ is ``approximately semifinite" (defined below).  Although neither of these properties is fully understood, we show that they are enjoyed by all semifinite algebras and hyperfinite algebras (with no summand of type $\text{I}_2$), plus others.  Thus the classification is stronger than Yeadon's theorem [Y1] for semifinite algebras (and its recent improvement in [JRS]), and the proof uses independent techniques.  Related to this, we examine the modular theory for positive projections from a von Neumann algebra onto a Jordan image of another von Neumann algebra, and use such projections to construct new $L^p$ isometries by interpolation.  Some complementary results and questions are also presented.
\end{abstract}

\maketitle

\section{Introduction}

In any class of Banach spaces, it is natural to ask about the isometries.  (Here an isometry is always assumed to be linear, but \textit{not} assumed to be surjective.)  $L^p$ function spaces are an obvious example, and their isometries have been understood for half a century.  To the operator algebraist, these classical $L^p$ spaces arise from commutative von Neumann algebras, and one may as well ask about isometries in the larger class of noncommutative $L^p$ spaces.  This question was considered by a number of authors, with a variety of assumptions; a succinct answer for semifinite algebras was given in 1981 by Yeadon [Y1].  In the recent paper [JRS], Yeadon's result was extended to the case where only the initial algebra is assumed semifinite.  We will call this the Generalized Yeadon Theorem (GYT) (Theorem \ref{T:gyt} below), as it was first proved by a simple modification of Yeadon's original argument.  But in its most general form, the classification of isometries between noncommutative $L^p$ spaces is still an open question.

Let us agree that ``$L^p$ isometry" will mean an isometric map $T: L^p(\mathcal{M}_1) \to L^p(\mathcal{M}_2)$, $1 \le p < \infty$, $p \ne 2$.  Adapting Watanabe's terminology ([W2], [W3]), we say that an $L^p$ isometry is \textbf{typical} if there are
\begin{enumerate}
\item a normal Jordan *-monomorphism $J: \mathcal{M}_1 \to \mathcal{M}_2$, \item a partial isometry $w \in \mathcal{M}_2$ with $w^*w = J(1)$, and 
\item a (not necessarily faithful) normal positive projection $P: \mathcal{M}_2 \to J(\mathcal{M}_1)$, 
\end{enumerate}
such that
\begin{equation} \label{E:typical}
T(\varphi^{1/p}) = w(\varphi \circ J^{-1} \circ P)^{1/p}, \qquad \forall\varphi \in (\mathcal{M}_1)_*^+.
\end{equation}
(Writing the projection as $P: \mathcal{M}_2 \to J(\mathcal{M}_1)$ will always mean that $P$ fixes $J(\mathcal{M}_1)$ pointwise.)  Here $\varphi^{1/p}$ is the generic positive element of $L^p(\mathcal{M}_1)$; see below for explanation.  Since any $L^p$ element is a linear combination of four positive ones, this completely determines $T$.  We will see in Section \ref{S:background} that typical isometries follow Banach's original classification paradigm for (classical) $L^p$ isometries, and may naturally be considered ``noncommutative weighted composition operators".

\begin{question} \label{T:conj}
Is every $L^p$ isometry typical?
\end{question}

Results in the literature offer evidence for an affirmative answer.  GYT and the structure theorem for $L^1$ isometries imply that an $L^p$ isometry with semifinite domain must be typical.  In [S1] typicality was proved for surjective $L^p$ isometries.  And the paper [JRS] shows typicality for $L^p$ isometries which are 2-isometries at the operator space level, with $J$ actually a homomorphism and $P$ a conditional expectation.

Our strategy here is the following.  First we use a theorem of Bunce and Wright [BW1] to show that $L^1$ isometries are typical.  (This has also been noted by Kirchberg [Ki].)  Then given an $L^p$ isometry, we try to form an associated $L^1$ isometry, apply typicality there, and deduce typicality for the original map.  It is not clear whether this procedure can work in general; it requires that continuous homogeneous positive bounded functions on $L^p(\mathcal{M}_1)_+$ which are additive on orthogonal elements must in fact be additive.  This is a natural variant of the extension property (EP) introduced by Watanabe [W2], so we call it EP$p$.

We will call a von Neumann algebra ``approximately semifinite" (AS) if it can be paved out by a net of semifinite subalgebras (see Section \ref{S:ep} for the precise definition).  The main results of Sections \ref{S:pone} through \ref{S:ep} are summarized in

\begin{theorem} \label{T:main} ${}$
\begin{enumerate}
\item All $L^1$ isometries are typical.
\item An $L^p$ isometry must be typical if $\mathcal{M}_1$ has EP$p$ and EP$1$; for positive $L^p$ isometries EP$1$ is sufficient.
\item An AS algebra with no summand of type $\text{I}_2$ has EP$p$ for any $p \in [1,\infty)$.
\item The class of AS algebras includes all semifinite algebras, hyperfinite algebras, factors of type $\text{III}_0$ with separable predual, and others.
\end{enumerate}
\end{theorem}

This is a stronger result than GYT, and the proofs are independent of Yeadon's paper.  At this time we do not know any non-AS algebra, or any factor other than $M_2$ which does not have EP$p$.  Further insight into these properties may help to resolve Question \ref{T:conj}, and they seem to merit investigation in their own right.

The converse of Question \ref{T:conj} is also interesting.

\begin{question} \label{T:conj2}
For a given normal Jordan *-monomorphism $J:\mathcal{M}_1 \to \mathcal{M}_2$ and normal positive projection $P:\mathcal{M}_2 \to J(\mathcal{M}_1)$, does \eqref{E:typical} extend linearly to an $L^p$ isometry?
\end{question}

If $J(\mathcal{M}_1)$ is a von Neumann algebra, then $P$ is a conditional expectation and the answer to Question \ref{T:conj2} is yes.  This is not entirely new, but we collect the necessary arguments in Section \ref{S:conditioned}.  Then in Section \ref{S:jordan} we compare the general situation.  We are able to construct new $L^p$ isometries from these data by interpolation, but now they seem to depend on the choice of a reference state.  However, in Section \ref{S:factor} we show that the dependence is removed, and Question \ref{T:conj2} can again be answered affirmatively, if $P$ factors through a conditional expectation from $\mathcal{M}_2$ onto $J(\mathcal{M}_1)''$.  Exactly this issue was addressed in a recent work of Haagerup and St\o rmer [HS2], although in a more general setting.  We extend their investigation in our specific case, guaranteeing the necessary factorization whenever $\mathcal{M}_1$ is AS.  Combining this with Theorem \ref{T:main}, we acquire a complete parameterization of the isometries from $L^p(\mathcal{M}_1)$ to $L^p(\mathcal{M}_2)$ whenever $\mathcal{M}_1$ is AS.

\bigskip

EP was proposed as a tool for $L^p$ isometries by Keiichi Watanabe, and I thank him for making his preprints available to me.  With his permission, a few of his unpublished results are incorporated here into Theorem \ref{T:semifinite} (and clearly attributed to him).  I am also grateful to Marius Junge and Zhong-Jin Ruan for helpful conversations and for showing me Theorem \ref{T:jrx} from [JRX].

\section{History and background} \label{S:background}

It is not plausible to review the theory of noncommutative $L^p$ spaces at length here.  The reader may want to consult [Te1], [N], [K1], [Ya] for details of the constructions mentioned below; [PX] also includes an extensive bibliography.  Our interest, aside from refreshing the reader's memory, lies largely in setting up convenient notation and explaining why ``typical" isometries are a natural generalization of previous results going back to the origins of the subject.

In fact the fundamental 1932 book of Banach ([B], IX.5) already listed the surjective isometries of $\ell^p$ and $L^p[0,1]$ (with respect to Lebesgue measure).  In the second case, an $L^p$ isometry $T$ is uniquely decomposed as a weighted composition operator:
\begin{equation} \label{E:wcomp}
T(f) = h \cdot (f \circ \varphi) = (\text{sgn} \: h) \cdot |h| \cdot (f \circ \varphi),
\end{equation}
where $h$ is a measurable function and $\varphi$ is a measurable (a.e.) bijection of $[0,1]$.  Clearly $|h|$ is related to the Radon-Nikod\'ym derivative for the change of measure induced by $\varphi$.  Although Banach did not prove this classification, he did make the key observation that isometries on $L^p$ spaces must preserve disjointness of support; i.e.
\begin{equation} \label{E:disj}
fg = 0 \iff T(f)T(g) = 0.
\end{equation}
We will see that the equations \eqref{E:wcomp} and \eqref{E:disj} provide a model for all succeeding classifications.

The extension to non-surjective isometries on general (classical) $L^p$ spaces was made in 1958 by Lamperti [L].  His description was similar, but he noted that generally the bijection $\varphi$ must be replaced by a ``composition" induced by a set-valued mapping, called a \textit{regular set isomorphism}.  See [L] or [FJ] for details; for many measure spaces [HvN] one can indeed find a (presumably more basic) point mapping.  As Lamperti pointed out, \eqref{E:disj} follows from a characterization of equality in the Clarkson inequality.  That is,
\begin{equation} \label{E:commclark}
\|f+g\|^p_p + \|f-g\|^p_p = 2(\|f\|^p_p + \|g\|^p_p) \iff fg=0.
\end{equation}
This method also works for some other function spaces ([L], [FJ]).

It is interesting to note how much of the analogous noncommutative machinery was in place at this time.  First observe that from the operator algebraic point of view a regular set isomorphism is more welcome than a point mapping, being a map on the projections in the associated $L^\infty$ algebra.  In terms of this (von Neumann) algebra, equation \eqref{E:disj} tells us that the underlying map between projection lattices preserves orthogonality.  Dye [D] had studied exactly such maps in the noncommutative setting a few years before, showing that they give rise to normal Jordan *-isomorphisms.  And Kadison's classic paper [Ka] had demonstrated the correspondence between normal Jordan *-isomorphisms and isometries.  Noncommutative $L^p$ spaces were around, too, but the isometric theory would wait for noncommutative formulations of \eqref{E:commclark}. 

Let us recall the definition of the noncommutative $L^p$ space ($1 \le p < \infty$) associated to a semifinite algebra $\mathcal{M}$ equipped with a given faithful normal semifinite tracial weight $\tau$  (simply called a ``trace" from here on).  The earliest construction seems to be due to Segal [Se].  Consider the set
$$\{T \in \mathcal{M} \mid \|T\|_p \triangleq \tau(|T|^p)^{1/p} < \infty \}.$$
It can be shown that $\| \cdot \|_p$ defines a norm on this set, so the completion is a Banach space, denoted $L^p(\mathcal{M}, \tau)$.  It also turns out that one can identify elements of the completion with unbounded operators; to be specific, all the spaces $L^p(\mathcal{M}, \tau)$ are subsets of the *-algebra of $\tau$-measurable operators $\mathfrak{M}(\mathcal{M}, \tau)$ [N].  Clearly $\tau$ is playing the role of integration.

Before stating Yeadon's fundamental classification for isometries of semifinite $L^p$ spaces, we recall that a \textit{Jordan} map on a von Neumann algebra is a *-linear map which preserves the Jordan product $x \bullet y \triangleq (1/2)(xy + yx)$.  (We denote this by $\bullet$ instead of $\circ$ since we use the latter for composition very frequently.)  The unfamiliar reader may be comforted to know that a normal Jordan *-monomorphism from one von Neumann algebra into another is the direct sum of a *-isomorphism $\pi$ and a *-antiisomorphism $\pi'$, where $s(\pi) + s(\pi') \ge 1$ [St1].  (We use $s$ and its variants $s_\ell, s_r$ for ``(left/right) support of" throughout the paper.)  This is frequently misinterpreted in the literature.  Part (but not all) of the confusion comes from the fact that the image is typically not multiplicatively closed; the simplest example is
$$J:M_2 \to M_4, \qquad x \mapsto \left( \begin{smallmatrix} x & 0 \\ 0 & x^t \end{smallmatrix} \right),$$
where $t$ is the transpose map.  Accordingly, we will refer to a Jordan image of a von Neumann algebra as a ``Jordan algebra" in order to remind the reader that it is closed under the Jordan product and not the usual product.  This slightly abusive terminology should cause no confusion; we will not need the abstract definitions of Jordan algebras, JW-algebras, etc. 

\begin{theorem} \label{T:yeadon} ([Y1], Theorem 2)
A linear map $T:L^p(\mathcal{M}_1, \tau_1) \to L^p(\mathcal{M}_2, \tau_2)$
is isometric if and only if there exists
\begin{enumerate}
\item a normal Jordan *-monomorphism $J: \mathcal{M}_1 \to \mathcal{M}_2$,
\item a partial isometry $w \in \mathcal{M}_2$ with $w^*w = J(1)$, and
\item a positive self-adjoint operator $B$ affiliated with $\mathcal{M}_2$ such that the spectral projections of $B$ commute with $J(\mathcal{M}_1)$, $s(B) = J(1)$, and $\tau_1(x) = \tau_2(B^p J(x))$ for all $x\in \mathcal{M}_1^+$,
\end{enumerate}
all satisying
\begin{equation} \label{E:ythm}
T(x) = w B J(x), \qquad \forall x \in \mathcal{M}_1 \cap L^p(\mathcal{M}_1, \tau_1)
\end{equation}
Moreover, $J,B,$ and $w$ are uniquely determined by $T$.
\end{theorem}

Note the striking resemblance between \eqref{E:wcomp} and \eqref{E:ythm} - again, $B$ is related to a (noncommutative) Radon-Nikod\'ym derivative.

\bigskip

Traciality is essential in the construction of $L^p(\mathcal{M}, \tau)$, so another method is required for the general case.  We proceed by analogy: if $\mathcal{M}$ is supposed to be a noncommutative $L^\infty$ space, the associated $L^1$ space should be the predual $\mathcal{M}_*$.  This is not given as a space of operators, so it is not clear where the $p$th roots are.  Later, in Section \ref{S:conditioned}, we will discuss Kosaki's interpolation method [K1].  Here we recall the first construction, due to Haagerup ([H1],[Te1]), which goes as follows.  Choose a faithful normal semifinite weight $\varphi$ on $\mathcal{M}$.  The crossed product $\widetilde{\mathcal{M}} \triangleq \mathcal{M} \rtimes_{\sigma^\varphi} \mathbb{R}$ is semifinite, with canonical trace $\bar{\tau}$ and trace-scaling dual action $\theta$.  Then $\mathcal{M}_*$ can be identified, as an ordered vector space and as an $\mathcal{M}-\mathcal{M}$ bimodule, with the $\bar{\tau}$-measurable operators $T$ affiliated with $\widetilde{\mathcal{M}}$ satisfying $\theta_s(T) = e^{-s}T$.  We may simply transfer the norm to this set of operators, and we denote this space by $L^1(\mathcal{M})$.  Of course, because of the identification with $\mathcal{M}_*$, $L^1(\mathcal{M})$ does not depend (up to isometric isomorphism) on the choice of $\varphi$.

We will use the following intuitive notation: for $\psi \in \mathcal{M}_*^+$, we also denote by $\psi$ the corresponding operator in $L^1(\mathcal{M})_+$.  (In the original papers this was written as $h_\psi$, but several other notations are in use -- it is $\Delta_{\psi, \varphi} \otimes \lambda$ in the crossed product construction of the last paragraph.  Some advantages and applications of our convention, called a \textit{modular algebra}, are demonstrated in ([C3], V.B.$\alpha$), [Ya], [FT], [JS], [S2].)  Recall that $x \in \mathcal{M}$ and $\psi \in \mathcal{M}_*$ are said to \textit{commute} if the functionals $x\psi = \psi(\cdot x)$ and $\psi x = \psi(x \cdot)$ agree.  Of course this is nothing but the requirement that $x \in \mathcal{M} \subset \widetilde{\mathcal{M}}$ and $\psi \in L^1(\mathcal{M})$ commute as operators.  We also have the useful relations
$$\varphi^{it} \psi^{-it} = (D\varphi: D \psi)_t; \qquad \psi^{it} x \psi^{-it} = \sigma^\psi_t(x)$$
for $\varphi, \psi \in \mathcal{M}_*^+, \: s(\varphi) \le s(\psi), \: x \in s(\psi) \mathcal{M} s(\psi).$  

Now we set $L^p(\mathcal{M})$ ($1 \le p < \infty$) to be the set of $\bar{\tau}$-measurable operators $T$ for which $\theta_s(T) = e^{-s/p}T$, and defining a norm $\|T\|_p = \||T|^p\|_1^{1/p}$ gives us a Banach space.  As a space of operators, $L^p(\mathcal{M})$ is still ordered, and any element is a linear combination of four positive ones.  This all agrees with our previous construction in case $\mathcal{M}$ is semifinite: the identification is
\begin{equation} \label{E:twolp}
L^p(\mathcal{M}, \tau)_+ \ni h \leftrightarrow (\tau_{h^p})^{1/p} \in L^p(\mathcal{M})_+.
\end{equation}
(Here $\tau_h(x) \triangleq \tau(hx)$; more generally $\varphi_h(x) \triangleq \varphi(hx) = \varphi(xh)$ whenever $\varphi$ and $h$ commute.)  But the reader should appreciate the paradigm shift: now $L^1$ elements are ``noncommutative measures".  Any theory for functorially producing $L^p$ spaces from von Neumann algebras (i.e. without arbitrarily choosing a base measure) is forced into such a construction, as von Neumann algebras do not come with distinguished measures unless the algebra is a direct sum of type I or $\text{II}_1$ factors.  It is more correct to think of a von Neumann algebra as determining a measure \textit{class} (in the sense of absolutely continuity), and this generates an $L^p$ space of measures directly.  See [S2] for more discussion.

\bigskip

Now we revisit Theorem \ref{T:yeadon}.  The operator $B$ commutes with $J(\mathcal{M}_1)$, and so when $\mathcal{M}_1$ is finite, the linear functional $\varphi \triangleq |T(\tau_1^{1/p})|^p = (\tau_2)_{B^p}$ commutes with $J(\mathcal{M}_1)$.  (Equivalently, the restriction of $\varphi$ to $J(\mathcal{M}_1)''$ is a finite trace.)  Formulated in this way, Theorem \ref{T:yeadon} extends to the case where $\mathcal{M}_2$ is not assumed semifinite (and $\mathcal{M}_1$ is not assumed finite).  The result, which we call GYT for ``Generalized Yeadon Theorem", was first noted in [JRS] but will be reproven as Theorem \ref{T:cgyt}.

\begin{theorem} \label{T:gyt} ([JRS])
Let $\mathcal{M}_1, \mathcal{M}_2$ be von Neumann algebras, $\tau$ a fixed trace on $\mathcal{M}_1$, and $1 \le p < \infty, \; p \ne 2$.  If $T: L^p(\mathcal{M}_1, \tau) \to L^p(\mathcal{M}_2)$ is an isometry, then there are, uniquely,
\begin{enumerate}
\item a normal Jordan *-monomorphism $J: \mathcal{M}_1 \to \mathcal{M}_2$,
\item a partial isometry $w \in \mathcal{M}_2$ with $w^*w = J(1)$, and
\item a normal semifinite weight $\varphi$ on $\mathcal{M}_2$, which commutes with $J(\mathcal{M}_1)''$ and satisfies $s(\varphi) = J(1)$, $\varphi(J(x)) = \tau(x)$ for all $x \in (\mathcal{M}_1)_+$,
\end{enumerate}
all satisfying
\begin{equation} \label{E:gyt}
T(x) = w\varphi^{1/p} J(x),\: \forall x \in \mathcal{M}_1 \cap L^p(\mathcal{M}_1, \tau).
\end{equation}
\end{theorem}

\textit{Remark.} An operator interpretation of $\tau$ and $\varphi$ requires a little more explanation when they are unbounded functionals ([Ya],[S2]), or one can rewrite \eqref{E:gyt} as
\begin{equation} \label{E:gyt2}
T(h^{1/p}) = w(\varphi_{J(h)})^{1/p}, \qquad h \in \mathcal{M}_1 \cap L^1(\mathcal{M}_1, \tau)_+,
\end{equation}
and extend by linearity.  We will use \eqref{E:gyt2} in the sequel.

\bigskip

For Theorems \ref{T:yeadon} and \ref{T:gyt}, a key ingredient of the proofs is the equality condition in the Clarkson inequality for \textit{noncommutative} $L^p$ spaces.  Yeadon [Y1] showed this for semifinite von Neumann algebras; a few years later Kosaki [K2] proved it for arbitrary von Neumann algebras with $2<p<\infty$; and only recently Raynaud and Xu [RX] obtained a general version (relying on Kosaki's work).  It plays a role in this paper as well.

\begin{theorem} \label{T:clarkson} [RX] (Equality condition for noncommutative Clarkson inequality)

For $\xi, \eta \in L^p(\mathcal{M})$, $0 < p < \infty$, $p \ne 2$,
\begin{equation} \label{E:clarkson}
\|\xi + \eta\|^p + \|\xi - \eta\|^p = 2(\| \xi \|^p + \| \eta \|^p) \iff \xi \eta^* = \xi^* \eta = 0.
\end{equation}
\end{theorem}

We remind the reader that $L^p$ elements have left and right support projections in $\mathcal{M}$.  Since $s_\ell(\xi) \perp s_\ell(\eta) \iff \xi^* \eta = 0$ as elements of Haagerup's $L^p$ space, we will call pairs satisfying the conditions of \eqref{E:clarkson} \textbf{orthogonal}.  (This can be interpreted in terms of $L^{p/2}$-valued inner products, see [JS].)  We mentioned earlier that isometries of classical $L^p$ spaces preserve disjointness of support: Theorem \ref{T:clarkson} tells us that all $L^p$ isometries actually preserve orthogonality, which is disjointness of left \textit{and} right supports.

Comparing \eqref{E:typical} and \eqref{E:wcomp}, one sees that typical $L^p$ isometries correspond to a noncommutative interpretation of Banach's classification result for $L^p(0,1)$.  One may think of the partial isometry $w$ as a ``noncommutative function of unit modulus" (corresponding to sgn $h$), and the precomposition with $J^{-1} \circ P$ as a `` noncommutative isometric composition operator" (corresponding to $f \mapsto |h| \cdot (f \circ \varphi)$).

\section{$L^1$ isometries} \label{S:pone}

The starting point for our investigation is the following (paraphrased) result of Bunce and Wright.  Recall that an \textit{o.d. homomorphism} $(\mathcal{M}_1)_* \to (\mathcal{M}_2)_*$ is a linear homomorphism which is positive and preserves orthogonality between positive functionals.

\begin{theorem} \label{T:bw} ( $\sim$ [BW1], Theorem 2.6)
If $T: (\mathcal{M}_1)_* \to (\mathcal{M}_2)_*$ is an o.d. homomorphism, then the map
$$J: s(\varphi) \to s(T(\varphi))$$
is well-defined and extends to a normal Jordan *-homomorphism.  We have $T^*(1)$ central in $\mathcal{M}_1$, and 
\begin{equation} \label{E:bw}
T(\varphi)(J(x)) = \varphi(T^*(1) x), \qquad \forall x \in \mathcal{M}_1, \: \forall \varphi \in (\mathcal{M}_1)_*.
\end{equation}
\end{theorem}

Consider the case where $T$ is a \textit{positive} isometry of $(\mathcal{M}_1)_*$ into $(\mathcal{M}_2)_*$.  Then $T$ is an o.d. homomorphism, as the equality condition of the Clarkson inequality shows:
\begin{align}
\notag \varphi \perp \psi &\Rightarrow \|\varphi + \psi\| + \|\varphi - \psi\| = 2(\|\varphi\| + \|\psi\|) \\
\notag &\Rightarrow \|T(\varphi) + T(\psi)\| + \|T(\varphi) - T(\psi)\| = 2(\|T(\varphi)\| + \|T(\psi)\|) \\
\notag &\Rightarrow T(\varphi) \perp T(\psi).
\end{align}
Applying Theorem \ref{T:bw} and equation \eqref{E:bw}, we first note that
\begin{equation} \label{E:unital}
\|\varphi\| = \|T(\varphi)\| = T(\varphi)(s(T(\varphi))) = T(\varphi)(J(s(\varphi))) = \varphi(T^*(1) s(\varphi))
\end{equation}
for each $\varphi \in (\mathcal{M}_1)_*$, which is only possible if $J$ is a monomorphism and $T^*(1) =1.$  By \eqref{E:bw} we have that $T^* \circ J = \text{id}_{\mathcal{M}_1}$.  Then $P \triangleq J \circ T^*$ is a normal positive projection from $\mathcal{M}_2$ onto $J(\mathcal{M}_1)$, $T^* = J^{-1} \circ P$, and $T$ is $(J^{-1}\circ P)_*$.

The following observation will be useful.  Since $P = J \circ T^*$, the supports of $P$ and $T^*$ are the same.  But $s(T^*)$ is the smallest projection in $\mathcal{M}_2$ such that for all $x \in (\mathcal{M}_2)_+, \: \varphi \in (\mathcal{M}_1)_*^+,$
$$T(\varphi)(x) = \varphi(T^*(x)) = \varphi(T^*(s(T^*)xs(T^*)) = T(\varphi)(s(T^*)x s(T^*)).$$
Thus
\begin{equation} \label{E:suppofp}
s(P) = s(T^*) = \sup_\varphi s(T(\varphi)) = \sup_\varphi J(s(\varphi)) = J(1) = P(1).
\end{equation}

Now consider an isometry $T$ from $L^1(\mathcal{M}_1)$ to $L^1(\mathcal{M}_2)$ which is not necessarily positive.  Let $\varphi, \psi \in (\mathcal{M}_1)_*^+$ be arbitrary, and let the polar decompositions be
$$T(\varphi) = u|T(\varphi)|; \quad T(\psi) = v|T(\psi)|; \quad T(\varphi+ \psi) = w|T(\varphi + \psi)|.$$
So $u^*u = s(\varphi)$, and similarly for the others.  Then
$$|T(\varphi + \psi)| = w^*T(\varphi + \psi) = w^*(T(\varphi) + T(\psi)) = w^*u|T(\varphi)| + w^*v|T(\psi)|.$$
View both sides as linear functionals and evaluate at 1:
$$|T(\varphi + \psi)|(1) = |T(\varphi)|(w^*u) + |T(\psi)|(w^*v) \le |T(\varphi)|(u^*u) + |T(\psi)|(v^*v)$$ $$ = \|T(\varphi)\| + \|T(\psi)\| = \|\varphi\| + \|\psi\| = \|\varphi + \psi\| = \|T(\varphi + \psi)\| = |T(\varphi + \psi)|(1).$$
Apparently the inequality is an equality, which implies by Cauchy-Schwarz that $w \, s(|T(\varphi)|)= u, \: w \, s(|T(\psi)|)= v.$  It follows that any partial isometry occurring in the polar decomposition of some $T(\varphi)$ is a reduction of a largest partial isometry $w$, with $s_r(w) = \vee \{s_r(T(\varphi)) \mid \varphi \in (\mathcal{M}_1)_*^+ \}$, so that $T(\varphi) = w|T(\varphi)|$ for any $\varphi \in (\mathcal{M}_1)^+_*$.  Then $\xi \mapsto w^*T(\xi)$ is a positive linear isometry, and we may use the previous argument to obtain the decomposition $T(\xi)=w (J^{-1} \circ P)_*(\xi).$  Necessarily by \eqref{E:suppofp} $w^*w = J(1)$; if there is a faithful normal state $\rho$ on $\mathcal{M}_1$, $w$ occurs in the polar decomposition of $T(\rho)$.  

We have shown that

\begin{theorem} \label{T:lone}
An $L^1$ isometry is typical.
\end{theorem}

A different proof of this can be found in Kirchberg ([Ki], Lemma 3.6).

\bigskip

\textit{Remark.} Using GYT and Theorem \ref{T:lone}, it is not difficult to prove that all $L^p$ isometries with semifinite domain are typical.  For if \eqref{E:gyt2} holds, one may define the positive $L^1$ isometry
$$T': \tau_{x} \mapsto \varphi_{J(x)}, \qquad x \in L^1(\mathcal{M}_1, \tau)_+$$
and deduce typicality for $T$ from that of $T'$.  As the development of this paper is intended to be independent of GYT, we derive this fact (and GYT) formally in Section \ref{S:ep}.

\section{$L^p$ isometries, $p >1$} \label{S:p>1}

Now consider an $L^p$ isometry $T$ with $1 < p < \infty$, $p \ne 2$.  Define
$$\bar{T}: L^p(\mathcal{M}_1)_+ \to L^p(\mathcal{M}_2)_+, \qquad \bar{T}(\varphi^{1/p}) = |T(\varphi^{1/p})|.$$
Is this map linear on $L^p(\mathcal{M}_1)_+$?  To attack this question, we make the following

\begin{definition} \label{D:cfm}
A \textbf{continuous finite measure} (c.f.m.) on $L^p(\mathcal{M})_+$ ($1 \le p < \infty$) is a nonnegative real-valued function $\rho$ which satisfies
\begin{enumerate}
\item $\rho(\lambda \varphi^{1/p}) = \lambda \rho(\varphi^{1/p})$,
\item $\varphi \perp \psi \Rightarrow \rho(\varphi^{1/p} + \psi^{1/p}) = \rho(\varphi^{1/p}) + \rho(\psi^{1/p})$,
\item $\rho(\varphi^{1/p}) \le C\|\varphi^{1/p}\|$ for some $C<\infty$ (denote by $\|\rho\|$ the least such $C$),
\item $\varphi^{1/p}_n \to \varphi^{1/p} \Rightarrow \rho(\varphi^{1/p}_n) \to \rho(\varphi^{1/p})$,
\end{enumerate}
for $\psi, \varphi, \varphi_n \in \mathcal{M}^+_*, \: \lambda \in \mathbb{R}_+.$

A von Neumann algebra $\mathcal{M}$ will be said to have \textbf{EP$p$} (extension property for $p$) if every c.f.m. $\rho$ on $L^p(\mathcal{M})_+$ is additive.  This implies that $\rho$ extends uniquely to a continuous linear functional on all of $L^p(\mathcal{M})$ and thus may be identified with an element of $L^q(\mathcal{M})_+$ ($1/p + 1/q = 1$).
\end{definition}

\textit{Remark.}  These definitions are adapted from [W2], where c.f.m. are defined on $L^1(\mathcal{M})_+$ only (and $C = 1$, which is inconsequential).  Thus Watanabe's EP corresponds to EP1 in our context.

\bigskip

Returning to $\bar{T}$, we see that each element $\psi^{1/q} \in L^p(\mathcal{M}_2)_+$ generates a c.f.m. on $L^p(\mathcal{M}_1)_+$ via $\varphi^{1/p} \mapsto <\bar{T}(\varphi^{1/p}), \psi^{1/q}>.$  The only nontrivial conditions to check are (2) and (4).  $\bar{T}$ preserves orthogonality, so it is additive on orthogonal elements, proving (2).  (4) follows from a result of Raynaud ([R], Lemma 3.2) on the continuity of the absolute value map in $L^p$, $0<p<\infty$.  The same lemma also shows that the map
\begin{equation} \label{E:raynaud}
L^p_+ \to L^q_+, \qquad \varphi^{1/p} \mapsto \varphi^{1/q} \qquad (0<p,q<\infty)
\end{equation}
is continuous, which will be useful shortly.

If $\mathcal{M}_1$ has EP$p$, then the c.f.m. generated by $\psi^{1/q}$ must be evaluation at some positive element of $L^q(\mathcal{M}_1)$.  We denote this element by $\pi(\psi^{1/q})$, so
\begin{equation} \label{E:pi}
<\varphi^{1/p}, \pi(\psi^{1/q})> = <\bar{T}(\varphi^{1/p}), \psi^{1/q}>.
\end{equation}
Now for all $\varphi^{1/p} \in L^p(\mathcal{M}_1)_+$,
$$<\varphi^{1/p}, \pi(\psi^{1/q})> = <|T(\varphi^{1/p})|, \psi^{1/q}> \le \|T(\varphi^{1/p})\| \|\psi^{1/q}\| = \|\varphi^{1/p}\| \|\psi^{1/q}\|,$$
so $\pi$ is \textit{norm-decreasing}.  And
\begin{align}
\notag <\varphi^{1/p}, \pi(\psi_1^{1/q} + \psi_2^{1/q})> &= <\bar{T}(\varphi^{1/p}), \psi_1^{1/q} + \psi_2^{1/q}> \\
\notag &= <\bar{T}(\varphi^{1/p}), \psi_1^{1/q}> + <\bar{T}(\varphi^{1/p}), \psi_2^{1/q}> \\
\notag &= <\varphi^{1/p}, \pi(\psi_1^{1/q})> + <\varphi^{1/p}, \pi(\psi_2^{1/q})>,
\end{align}
so $\pi$ is \textit{linear}.  Also denote by $\pi$ the unique linear extension to all of $L^q(\mathcal{M}_2)$.  Now by \eqref{E:pi}, $\bar{T}$ agrees with $\pi^*$ on $L^p(\mathcal{M}_1)_+$.  In particular, $\bar{T}$ is additive.

A symmetric argument shows that the map
\begin{equation} \label{E:reverse}
\varphi^{1/p} \mapsto |T(\varphi^{1/p})^*|, \qquad \varphi^{1/p} \in L^p(\mathcal{M}_1)_+
\end{equation}
is additive.  Knowing that these two maps are additive allows us to find one of the ingredients of typicality, the partial isometry.

Choose $\varphi, \psi \in \mathcal{M}_*^+$, and let the polar decompositions be
$$T(\varphi^{1/p}) = u|T(\varphi^{1/p})|; \quad T(\psi^{1/p}) = v|T(\psi^{1/p})|;$$ $$ T(\varphi^{1/p}+ \psi^{1/p}) = w|T(\varphi^{1/p} + \psi^{1/p})|.$$
We calculate
$$\left[u|T(\varphi^{1/p})|^{1/2} - w|T(\varphi^{1/p})|^{1/2}\right] \left[u|T(\varphi^{1/p})|^{1/2} - w|T(\varphi^{1/p})|^{1/2}\right]^*$$ $$ + \left[v|T(\psi^{1/p})|^{1/2} - w|T(\psi^{1/p})|^{1/2}\right]  \left[u|T(\psi^{1/p})|^{1/2} - w|T(\psi^{1/p})|^{1/2}\right]^*$$
$$= u|T(\varphi^{1/p})|u^* + w|T(\varphi^{1/p})|w^* - u|T(\varphi^{1/p})|w^* - w|T(\varphi^{1/p})|u^*$$ $$ + v|T(\psi^{1/p})|v^* + w|T(\psi^{1/p})|w^* - v|T(\psi^{1/p})|w^* - w|T(\psi^{1/p})|v^*.$$
Now we use
$$u|T(\varphi^{1/p})|u^* + v|T(\psi^{1/p})|v^* = w|T(\varphi^{1/p} + \psi^{1/p})|w^* = w|T(\varphi^{1/p})|w^* + w|T(\psi^{1/p})|w^*$$
(which follows from additivity of \eqref{E:reverse} and $\bar{T}$) on the first and fifth term, and
$$u|T(\varphi^{1/p})| + v|T(\psi^{1/p})| = w|T(\varphi^{1/p} + \psi^{1/p})| = w|T(\varphi^{1/p})| + w|T(\psi^{1/p})|$$
(which follows from additivity of $T$ and $\bar{T}$) on the third and seventh, and fourth and eighth.  This gives
$$w|T(\varphi^{1/p})|w^* + w|T(\varphi^{1/p})|w^* - w|T(\varphi^{1/p})|w^* -w|T(\varphi^{1/p})|w^*$$ $$+ w|T(\psi^{1/p})|w^* + w|T(\psi^{1/p})|w^* -w|T(\psi^{1/p})|w^* - w|T(\psi^{1/p})|w^* = 0.$$
We conclude that
$$ u|T(\varphi^{1/p})|^{1/2} = w|T(\varphi^{1/p})|^{1/2}, \qquad v|T(\psi^{1/p})|^{1/2} = w|T(\psi^{1/p})|^{1/2},$$
which means that $u$ and $v$ are restrictions of $w$.  Then there is a largest partial isometry $w$ with $T(\varphi^{1/p})= w|T(\varphi^{1/p})|$ for all $\varphi^{1/p} \in L^p(\mathcal{M}_1)_+$.  This means that the map
$$ \xi \mapsto w^*T(\xi), \qquad \xi \in L^p(\mathcal{M}_1) $$
is a \textit{positive} linear isometry.

So it suffices to show typicality for a positive $L^p$ isometry $T$.  We will now assume that $\mathcal{M}_1$ has EP1.  Consider the map
\begin{equation} \label{E:T'isom}
T': (\mathcal{M}_1)_*^+ \to (\mathcal{M}_2)_*^+; \qquad \varphi \to T(\varphi^{1/p})^p.
\end{equation}
By mimicking the argument given above for $\bar{T}$, we may use EP1 to show that $T'$ is additive.  (Each $h$ in $(\mathcal{M}_2)_+$ generates a c.f.m. on $(\mathcal{M}_1)_*^+$ by $\varphi \mapsto T'(\varphi)(h)$.  If we denote by $\pi'(h)$ the corresponding element of $(\mathcal{M}_1)_*^+$, then $\pi'$ is linear and extends to all of $\mathcal{M}_2$.  We have that $T'$ is the restriction of $(\pi')_*$ to $(\mathcal{M}_1)_*^+$.)

Now extend $T'$ linearly to all of $\mathcal{M}_*$ (as $(\pi')_*$).  Apparently $T'$ is an o.d. homomorphism (rememeber that $T$ preserves orthogonality), so we may apply Theorem \ref{T:bw}.  Since $T'$ is isometric on $(\mathcal{M}_1)_*^+$ by \eqref{E:T'isom}, \eqref{E:unital} again shows that $J$ is a monomorphism and $(T')^*(1) = 1$.  Then $(T')^* \circ J = \text{id}_{\mathcal{M}_1}$, and $P \triangleq J \circ (T')^*$ is a normal positive projection from $\mathcal{M}_2$ onto $J(\mathcal{M}_1)$.

Finally, notice that for any $x \in \mathcal{M}_2$,
$$T(\varphi^{1/p})^p(x)= T'(\varphi)(x) = \varphi ((T')^*(x)) = \varphi (J^{-1} \circ J \circ (T')^*(x)) = \varphi(J^{-1} \circ P(x)).$$
Therefore $T(\varphi^{1/p}) = (\varphi \circ J^{-1} \circ P)^{1/p}.$  We have shown

\begin{theorem} \label{T:ep}
Let $T$ be an isometry from $L^p(\mathcal{M}_1)$ to $L^p(\mathcal{M}_2)$ ($1<p<\infty$, $p \ne 2$), and assume $\mathcal{M}_1$ has EP$p$ and EP1.  Then $T$ is typical.  If $T$ is positive, then EP1 alone is sufficient to conclude typicality.
\end{theorem}

\section{EP$p$ algebras} \label{S:ep}

Probably the reader is already wondering: Which von Neumann algebras have EP$p$?  
This an $L^p$ version of an old question of Mackey on linear extensions of measures on projections.  The most relevant formulation is the following: suppose $\mu$ is a bounded nonnegative real-valued function on the projections in a von Neumann algebra $\mathcal{M}$ which is $\sigma$-additive on orthogonal projections.  Is $\mu$ the restriction of a normal linear functional?  The answer is yes, provided $\mathcal{M}$ has no summand of type $\text{I}_2$.  This was achieved in stages by Gleason [G], Christensen [Ch], Yeadon [Y2]; for a very general result see [BW2].  It is tempting to expect the same answer for EP$p$ - and this would resolve Question \ref{T:conj} affirmatively, by Theorem \ref{T:ep} - but there is no obvious $L^p$ analogue for the lattice of projections in a von Neumann algebra.  For example, a state with trivial centralizer [HT] cannot be written in any way as the sum of two orthogonal positive normal linear functionals.

At the other extreme, a trace $\tau$ allows us to embed the $\tau$-finite elements densely into the predual while preserving orthogonality.  This leads to Theorem \ref{T:semifinite}, due essentially to Watanabe ([W3], Lemma 6.7 and Theorem 6.9).  Working with EP1, he proved the first part for sequences and the second for finite von Neumann algebras.  With his permission, we incorporate his proof in the one given here.

We need a little preparation.

\begin{definition} \label{D:paving}
Let $\mathcal{M}$ be a von Neumann algebra and $\{E_\alpha\}$ a net of normal conditional expectations onto increasing subalgebras $\{\mathcal{M}_\alpha\}$.  We do not assume that the $E_\alpha$ are faithful, but we do require that $s(E_\alpha) = E_\alpha(1)$ (which is the unit of $\mathcal{M}_\alpha$).  Assume further that $\cup \mathcal{M}_\alpha$ is $\sigma$-weakly dense in $\mathcal{M}$, and $E_\alpha \circ E_\beta = E_\alpha$ for $\alpha < \beta$.  Then we say that $\mathcal{M}$ is \textbf{paved out} by $\{\mathcal{M}_\alpha, E_\alpha \}$.
\end{definition}

\begin{theorem} \label{T:tsukada} ([Ts], Theorem 2)  Let $\mathcal{M}$ be paved out by $\{\mathcal{M}_\alpha, E_\alpha \}$.
\begin{enumerate}
\item For any $\theta \in \mathcal{M}^+_*$, $(\theta \circ E_\alpha) \to \theta$ in norm.
\item For any $x \in \mathcal{M}$, $E_\alpha(x)$ converges strongly to $x$.  (We write $E_\alpha(x) \overset{s}{\to} x$.)
\end{enumerate}
\end{theorem}

Theorem \ref{T:tsukada} is proved by Tsukada ([Ts], Theorem 2) in a slightly different guise.  He does not start by assuming that $\cup \mathcal{M}_\alpha$ is $\sigma$-weakly dense in $\mathcal{M}$, but reduces to this case by requiring that all $E_\alpha$ preserve some faithful normal semifinite weight.  He also requires that the $E_\alpha$ are faithful, so let us show how his proof may be altered to handle the weaker assumption $s(E_\alpha) = E_\alpha(1)$.

The faithfulness is used in showing that if $E_\alpha(x) = 0$ for all $\alpha$, then $x =0$.  First Tsukada deduces that $E_{\alpha} (x^*x) = 0$ for any $\alpha$, and of course the faithfulness of a single $E_{\alpha}$ immediately implies that $x^*x=0.$  Without faithfulness, we obtain that
$$E_{\alpha}[s(E_{\alpha}) (x^* x) s(E_{\alpha})] = E_{\alpha}(x^*x) = 0 \Rightarrow s(E_{\alpha}) x^*x s(E_{\alpha}) = 0.$$
By assumption $s(E_\alpha) \nearrow 1,$ so we still conclude $x^*x = 0$.

\begin{theorem} \label{T:semifinite} Let $1 \le p < \infty$.
\begin{enumerate}
\item If $\mathcal{M}$ is paved out by $\{\mathcal{M}_\alpha, E_\alpha\}$, and each $\mathcal{M}_\alpha$ has EP$p$, then $\mathcal{M}$ has EP$p$.
\item A semifinite von Neumann algebra with no summand of type $\text{I}_2$ has EP$p$.
\end{enumerate}
\end{theorem}

\begin{proof}
Assume the hypotheses of (1) and let $\rho$ be a c.f.m. on $L^p(\mathcal{M})_+$.  Then
$$\rho_\alpha(\varphi^{1/p}) \triangleq \rho((\varphi \circ E_\alpha)^{1/p}), \qquad \varphi^{1/p} \in L^p(\mathcal{M}_\alpha)_+,$$
defines a c.f.m on $L^p(\mathcal{M}_\alpha)_+$.  (Note that $\rho_\alpha$ is continuous because the map
\begin{equation} \label{E:firste}
\varphi^{1/p} \mapsto (\varphi \circ E_\alpha)^{1/p}
\end{equation}
generates an isometric embedding $L^p(\mathcal{M}_\alpha) \hookrightarrow L^p(\mathcal{M})$, as reviewed in Section \ref{S:conditioned}.  Since $\varphi$ and $\varphi \circ E_\alpha$ have the same support in $\mathcal{M}_\alpha \subset \mathcal{M}$, $\rho_\alpha$ is additive on orthogonal elements.)  We have assumed that $\mathcal{M}_\alpha$ has EP$p$, so there is $\psi^{1/q}_\alpha \in L^q(\mathcal{M}_\alpha)_+$, $\|\psi^{1/q}_\alpha\| \le \|\rho\|$, with $\rho_\alpha(\varphi^{1/p}) = <\varphi^{1/p}, \psi^{1/q}_\alpha>$.  Now for any $\theta^{1/p} \in L^p(\mathcal{M})_+$, we have $(\theta \circ E_\alpha)^{1/p} \to \theta^{1/p}$ in norm.  This follows from Theorem \ref{T:tsukada}(1) and the continuity of \eqref{E:raynaud}.

We invoke the continuity of $\rho$ to calculate
$$\rho(\theta^{1/p})=\lim \rho( (\theta \circ E_\alpha)^{1/p}) = \lim \rho_\alpha((\theta \mid_{\mathcal{M}_\alpha} )^{1/p})$$ $$ = \lim < (\theta \mid_{\mathcal{M}_\alpha})^{1/p}, \psi^{1/q}_\alpha> = \lim < (\theta \circ E_\alpha)^{1/p}, (\psi_\alpha \circ E_\alpha)^{1/q}>.$$
The last equality depends on the fact that the family of inclusions \eqref{E:firste} also preserves duality, as mentioned in Section \ref{S:conditioned}.

These arguments show that
\begin{align}
\notag &<\theta^{1/p}, (\psi_\alpha \circ E_\alpha)^{1/q}> \\
\notag &{} \quad = <\theta^{1/p} - (\theta \circ E_\alpha)^{1/p}, (\psi_\alpha \circ E_\alpha)^{1/q}> + <(\theta \circ E_\alpha)^{1/p}, (\psi_\alpha \circ E_\alpha)^{1/q}> \\
\notag &{} \quad \to 0 + \rho(\theta^{1/p}).
\end{align}
(Note that $\|(\psi_\alpha \circ E_\alpha)^{1/q}\| = \|\psi_\alpha^{1/q}\|$ is bounded.)  Then $\rho$ is the limit of linear functionals and therefore linear itself, so $\mathcal{M}$ has EP$p$.

To prove part (2), first consider a \textit{finite} algebra $\mathcal{N}$ with normal faithful trace $\tau$ and no summand of type $\text{I}_2$.  Given a c.f.m. $\rho$, define the following measure on the projection lattice of $\mathcal{N}$: $\Phi(q) = \rho(q \tau^{1/p})$.  The continuity of $\rho$ implies that $\Phi$ is $\sigma$-additive.  By the result mentioned at the beginning of this section, there must be $\varphi \in \mathcal{M}_*^+$ with $\Phi(q) = \varphi(q)$.  Since $\mathcal{N}$ is finite, $\varphi$ is of the form $\tau_h$ for some $h \in L^1(\mathcal{N}, \tau)_+.$  We obtain
$$\rho(q\tau^{1/p}) = \tau_h(q).$$
Any element of $L^p(\mathcal{N}, \tau)_+$ is well-approximated by a finite positive linear combination of orthogonal projections, so $\rho$ being a c.f.m. gives us
$$\rho(k \tau^{1/p}) = \tau(hk), \qquad \forall k \in \mathcal{N}_+.$$
Now the map $k \tau^{1/p} \mapsto \tau(hk)$ is bounded (by $\|\rho\|$), so we must have $h \in L^q(\mathcal{N}, \tau)_+$.  That is,
$$\rho(k \tau^{1/p}) = < k \tau^{1/p}, h \tau^{1/q} >,$$
and so $\mathcal{N}$ has EP$p$.

Part (2) then follows from (1): given $\mathcal{M}$ semifinite, we may fix a faithful normal semifinite trace $\tau$ and notice that $\mathcal{M}$ is paved out by
$$\{q_\alpha \mathcal{M} q_\alpha, \: E_\alpha: x \mapsto q_\alpha x q_\alpha\},$$
where $q_\alpha$ runs over the lattice of $\tau$-finite projections.
\end{proof}

\textit{Remark 1.} Just as in Mackey's question, von Neumann algebras of type $\text{I}_2$ do not have EP$p$.  In $M_2$, for example, the manifold of one-dimensional projections is homeomorphic to the Riemann sphere $S^2$.  To extend (using Definition \ref{D:cfm}) to a c.f.m. on $L^p(M_2)_+$, a continuous nonnegative function $\rho$ on the sphere only needs to satisfy
$$\rho(p) + \rho(1-p) = \text{constant}, \qquad \forall p \in S^2.$$
(This is because 1 is the only element which may be written in more than one way as an orthogonal sum of positive elements.)  But typically such a c.f.m. will not be linear with respect to the vector space structure of $L^p(M_2)$.  The space of functions defined above is an infinite-dimensional real cone, but $L^q(M_2)_+$ has dimension four.

\textit{Remark 2.} If one wants to determine whether \textit{all} von Neumann algebras (without summands of type $\text{I}_2$) have EP$p$, it suffices to work only in the $\sigma$-finite case.  This follows from the fact that any von Neumann algebra $\mathcal{M}$ is paved out by 
$$\{p_\alpha \mathcal{M} p_\alpha, \: E_\alpha: x \mapsto p_\alpha x p_\alpha\},$$
where $p_\alpha$ runs over the lattice of $\sigma$-finite projections.

\bigskip

From Theorems \ref{T:semifinite}(2) and \ref{T:ep}, we see that an $L^p$ isometry with $\mathcal{M}_1$ semifinite (and lacking a type $\text{I}_2$ summand) must be typical.  After a preparatory lemma, we finally use this to give a new proof of GYT.

\begin{lemma} \label{T:stormer} 
Let $J: \mathcal{M}_1 \to \mathcal{M}_2$ be a normal Jordan *-monomorphism, $P: \mathcal{M}_2 \to J(\mathcal{M}_1)$ a normal positive projection, $x \in \mathcal{M}_1$, and $y \in \mathcal{M}_2$.
\begin{enumerate}
\item $\|P\| = 1.$
\item $P(J(x) \bullet y) = J(x) \bullet P(y).$
\item $P(J(x)yJ(x)) = J(x)P(y)J(x).$
\item $P(J(\mathcal{M}_1)' \cap \mathcal{M}_2) = J(\mathcal{Z}(\mathcal{M}_1))$.
\end{enumerate}
\end{lemma}

\begin{proof} Since $\|P(1)\|=1$, the first statement is a consequence of the corollary to the Russo-Dye Theorem [DR].  The next two statements are straightforward adaptations of ([St2], Lemma 4.1), but it will be useful to note here that the third follows from the second by the general Jordan algebra identity $aba = 2 a \bullet (a \bullet b) - a^2 \bullet b$.  The fourth is not new, but less explicit in our sources.  It follows from taking $z \in J(\mathcal{M}_1)' \cap \mathcal{M}_2$ and a projection $p \in \mathcal{M}_1$, and using the previous parts:
\begin{equation} \label{E:jcenter}
J(p) \bullet P(z) = P(J(p) \bullet z) = P(J(p)zJ(p)) = J(p)P(z)J(p).
\end{equation}
Applying $J^{-1}$ to \eqref{E:jcenter} and using the Jordan identity just mentioned gives
$$p \bullet [J^{-1} \circ P(z)] = p [J^{-1} \circ P(z)] p.$$
This implies that $J^{-1} \circ P(z) \in \mathcal{Z}(\mathcal{M}_1)$.
\end{proof}

Note that Lemma \ref{T:stormer}(2) is the Jordan version of the fact that conditional expectations are bimodule maps.

\begin{theorem} \label{T:cgyt}
Let $T$ be an $L^p$ isometry, and assume $(\mathcal{M}_1, \tau)$ is semifinite with no type $\text{I}_2$ summand.  Then GYT (Theorem \ref{T:gyt}) holds.
\end{theorem}

\begin{proof}
We first make the identification \eqref{E:twolp} between $L^p(\mathcal{M}_1, \tau)$ and $L^p(\mathcal{M}_1)$.  As just noted, $T$ is typical, so there are $w,J,P$ satisfying \eqref{E:typical}.  Letting $\varphi$ be the (necessarily normal and semifinite) weight $\tau \circ J^{-1} \circ P$, we have $\varphi(J(h)) = \tau(h)$ for $h \in (\mathcal{M}_1)_+$.  Equation \eqref{E:suppofp} guarantees that $s(\varphi) = J(1) = w^*w.$  It is left to show that $\varphi$ commutes with $J(\mathcal{M}_1)''$, to derive \eqref{E:gyt2}, and to show uniqueness of the data.

Because $\varphi$ may be unbounded, the commutation is more delicate than $\varphi J(x) = J(x) \varphi$.  The precise meaning is that $J(\mathcal{M}_1)'' \subset (\mathcal{M}_2)^\varphi$, the centralizer of $\varphi$; we need to show that $\sigma^\varphi$, which is defined on $s(\varphi)\mathcal{M}_2 s(\varphi)$, is the identity on $J(\mathcal{M}_1)''$.  One natural approach goes by Theorem \ref{T:hs}, but here we give a different argument.

Let $q$ be an arbitrary projection of $\mathcal{M}_1$, and let $s$ be the symmetry (=self-adjoint unitary) $1-2q$.  For $y \in (\mathcal{M}_2)_+$, we use Lemma \ref{T:stormer} to compute
$$\varphi(J(s)yJ(s)) = \tau \circ J^{-1} \circ P(J(s)yJ(s)) = \tau(s J^{-1} \circ P(y)s) = \varphi(y).$$
Thus $\varphi = \varphi \circ \text{Ad}\,J(s).$  By ([T2], Corollary VIII.1.4), for any $y \in J(1)\mathcal{M}_2J(1)$ and $t \in \mathbb{R}$,
\begin{align*}
\sigma^\varphi_t(y) = \sigma^{(\varphi \circ \text{Ad}\,J(s))}_t(y) &= \text{Ad}\,J(s) \circ \sigma^\varphi_t \circ \text{Ad}\,J(s)(y)\\
&= J(s) \sigma^\varphi_t(J(s)) \sigma^\varphi_t(y) \sigma^\varphi_t(J(s)) J(s).
\end{align*}
Since $y$ is arbitrary, we have that for each $t$, $J(s) \sigma^\varphi_t(J(s))$ belongs to the center of $J(1)\mathcal{M}_2 J(1)$.  Then
$$[J(s)\sigma^\varphi_t(J(s))]J(s) = J(s) [J(s)\sigma^\varphi_t(J(s))] = \sigma^\varphi_t(J(s))$$ $$ \Rightarrow J(s)\sigma^\varphi_t(J(s)) = \sigma^\varphi_t(J(s))J(s).$$
Central elements are fixed by modular automorphism groups, so
$$J(s)\sigma^\varphi_t(J(s)) = \sigma^\varphi_t(J(s))J(s) = \sigma_{-t}^\varphi[\sigma^\varphi_t(J(s))J(s)] = J(s)\sigma^\varphi_{-t}(J(s)).$$
Then
$$\sigma^\varphi_t(J(s)) = \sigma^\varphi_{-t}(J(s)) \Rightarrow \sigma^\varphi_{2t}(J(s)) = J(s).$$
So $\sigma^\varphi$ fixes all symmetries in $J(\mathcal{M}_1)$, so all projections in $J(\mathcal{M}_1)$, so all of $J(\mathcal{M}_1)$, and finally all of $J(\mathcal{M}_1)''$.  We will use this in the proof of Proposition \ref{T:asfactor}.

Now take any $h \in \mathcal{M}_1 \cap L^1(\mathcal{M}_1, \tau)_+, \: y \in \mathcal{M}_2$, and observe
$$\varphi_{J(h)}[y] = \tau \circ J^{-1} \circ P[J(h^{1/2}) y J(h^{1/2})] = \tau [h^{1/2} J^{-1} \circ P(y) h^{1/2}] = \tau_h \circ J^{-1} \circ P[y].$$
This implies 
$$w (\varphi_{J(h)})^{1/p} = w (\tau_h \circ J^{-1} \circ P)^{1/p} = T(h^{1/p}),$$
which is exactly \eqref{E:gyt2}.

For uniqueness, we repeat the argument from [JRS] for the convenience of the reader.  Suppose that we also have data $v, K, \psi$ verifying the hypotheses of GYT.  Then for any $\tau$-finite projection $q \in \mathcal{M}_1$, we have 
\begin{equation} \label{E:jk}
w (\varphi_{J(q)})^{1/p} = T(q) = v (\psi_{K(q)})^{1/p} \Rightarrow \varphi_{J(q)} = \psi_{K(q)}.
\end{equation}
Now take any projection $p \in \mathcal{M}_1$, and note that for all $\tau$-finite $q \le p$, 
$$\|\psi_{K(q)}\| = \|\varphi_{J(q)}\| = \varphi_{J(q)}(J(p)) = \psi_{K(q)}(J(p)).$$
This is only possible if $J(p) \ge K(q)$, and after taking the supremum over $q$ we get $J(p) \ge K(p)$.  A parallel argument gives the opposite inequality, implying $J=K$.  By \eqref{E:jk} we have $\psi_{J(q)} = \varphi_{J(q)}$ for all $\tau$-finite projections $q$, so $\varphi = \psi$ as weights.  That $w=v$ is now obvious.

\end{proof}

Of course GYT and typicality still hold on $\text{I}_2$ summands, by Yeadon's theorem and the remark at the end of Section \ref{S:pone}.  The uniqueness argument suggests the same statement for typical isometries, which we now prove.

\begin{proposition} \label{T:unique}
Any typical $L^p$ isometry can be written in the form \eqref{E:typical} for a unique triple $w,J,P$ satisfying $s(P) = P(1) (=J(1))$.
\end{proposition}

\begin{proof}  We always have that $s(P)$ commutes with $J(\mathcal{M}_1)$ (Lemma 1.2 of [ES]) and has the same central support as $P(1)$.  So if we consider the new Jordan *-monomorphism $J_0: x \mapsto J(x)s(P)$ and the new projection $P_0: y \mapsto P(y) s(P)$, we have $J^{-1} \circ P = J_0^{-1} \circ P_0$ and $s(P_0) = P_0(1)$.

To show uniqueness, suppose that an $L^p$ isometry can be written in terms of two triples $w,J,P$ and $w',J',P'$ satisfying all the necessary conditions.  By taking absolute values we get that $\varphi \circ J^{-1} \circ P = \varphi \circ J'^{-1} \circ P'$ for all $\varphi \in (\mathcal{M}_1)_*^+$, so we must have that the maps $J^{-1} \circ P$ and $J'^{-1} \circ P'$ agree.  Applying these to $J'(p),$ $p$ a projection in $\mathcal{M}_1$, gives
$$p = J^{-1} \circ P \circ J'(p), \qquad \text{or} \qquad J(p) = P \circ J'(p).$$
Using Lemma \ref{T:stormer}(3), we calculate
$$ P(J(p) J'(p) J(p)) = J(p) P (J'(p)) J(p) = J(p)J(p)J(p) = J(p) = P(J(p)).$$
But $J(p) J'(p) J(p) \le J(p)$, and $P$ is faithful on $J(1) \mathcal{M}_2 J(1)$.  Therefore the inequality is an equality, which implies $J'(p) \ge J(p)$.  The opposite inequality is derived symmetrically, so $J$ and $J'$ agree on projections and must agree everywhere.  Knowing this, it is easy to see that $P=P'$ and $w = w'$.
\end{proof}

Because of Proposition \ref{T:unique}, \textit{in the rest of the paper we will always assume that the support of a normal positive projection $P$ is equal to $P(1)$.}  This was incorporated into Definition \ref{D:paving} for the special case of conditional expectations, and we saw in \eqref{E:suppofp} that it already holds for all $P$ generated by Theorem \ref{T:bw}.

\bigskip

There are a few results in the literature which can be employed to establish EP$p$ in some type III von Neumann algebras.  Haagerup and St\o rmer ([HS1], Theorem 8.3) used a construction of Connes ([C1], Corollary 5.3.6) to show that factors of type $\text{III}_0$ with separable predual are paved out by $\text{II}_\infty$ algebras, so by Theorem \ref{T:semifinite} they all have EP$p$.  In another direction, we have the following theorem of Junge, Ruan, and Xu, which is a nontrivial modification of fundamental results for type III factors by Connes [C2] and Haagerup [H2].

\begin{theorem} \label{T:jrx} ([JRX], Theorem 4.3)
Let $\mathcal{M}$ be a hyperfinite type III von Neumann algebra.  Then there exist a normal faithful state $\varphi$ on $\mathcal{M}$ and an increasing sequence of $\varphi$-invariant normal faithful conditional expectations $\{E_k\}$ from $\mathcal{M}$ onto type I von Neumann subalgebras $\{\mathcal{N}_k\}$ of $\mathcal{M}$ such that $\cup \mathcal{N}_k$ is $\sigma$-weakly dense in $\mathcal{M}$.
\end{theorem}

This allows us to show

\begin{proposition} \label{T:hyperfinite}
Let $\mathcal{M}$ be a hyperfinite type III algebra.  Then $\mathcal{M}$ has EP$p$.
\end{proposition} 

\begin{proof}
This is a direct consequence of Theorems \ref{T:semifinite} and \ref{T:jrx}.  The only point which may not be obvious is that one can avoid $\text{I}_2$ summands.  That is easy to fix: given by Theorem \ref{T:jrx} the paving
$$E_k : \mathcal{M} \to \mathcal{N}_k$$
where $\mathcal{N}_k$ are type I, consider
$$E_k \otimes \text{id}: \mathcal{M} \otimes M_3 \to \mathcal{N}_k \otimes M_3.$$
Since $\mathcal{M}$ is isomorphic to $\mathcal{M} \otimes M_3$, $\{\mathcal{N}_k \otimes M_3\}$ can be identified with conditioned subalgebras of $\mathcal{M}$ which have no summands of type $\text{I}_2$.  Theorem \ref{T:semifinite} applies to the latter paving.
\end{proof}

\bigskip

\begin{definition}
A von Neumann algebra will be called \textbf{approximately semifinite (AS)} if it can be paved out by semifinite subalgebras.
\end{definition}
  
This terminology is an obvious analogy with the approximately finite-dimensional (AFD) algebras, but the author has been unable to find it in the literature.  In fact any AS algebra can be paved out by \textit{finite} subalgebras - since semifinite algebras can - but we wish to avoid the term ``approximately finite", which has another meaning.  So far we have seen that AS contains all semifinite algebras, hyperfinite algebras, and $\text{III}_0$ factors with separable predual.  It is closed under pavings, sums, and tensor products, so it also contains others: for example, the tensor product of a hyperfinite type III algebra with $L(\mathbb{F}_n)$, $n \ge 2$.  If we ignore $\text{I}_2$ summands, we have the inclusions of classes
$$\text{AS} \subset \text{EP}p \subset \text{vNa}.$$
As of this writing the author does not know if any of these inclusions are proper, or whether EP$p$ is independent of $p$ (we doubt the latter).  Regarding AS, we wish to point out that a paving may necessarily \textit{avoid} some subalgebras: the centralizer of an inner homogeneous state on a hyperfinite type III algebra is not properly contained in any other conditioned semifinite subalgebra (see [He], or 29.12 of [Str] for a related result).

At this point we have verified all the assertions in Theorem \ref{T:main}.  We should also mention that it is possible to show directly that $L^p$ isometries with AS domains are typical: pave out the domain with semifinite $L^p$ spaces, apply typicality to each subspace, and argue that the associated Jordan maps, partial isometries, and projections all converge in an appropriate sense.  Such a proof would presumably invoke Theorem \ref{T:hs}.

\section{$L^p$ isometries from *-(anti)isomorphisms and conditional expectations} \label{S:conditioned}

The projection $P$ occurring in our definition of typicality is formally very similar to a normal conditional expectation.  In this section we provide a construction of $L^p$ isometries associated to *-(anti)isomorphisms onto subalgebras which are the range of a normal conditional expectation.  Some of this material can be found in the literature, but it seems worthwhile to organize and justify the arguments.  (See Section 2 of [J], or Section 3 of [W1] for related discussions.)  In the succeeding sections we will try to formulate the analogous theory for $P$ and apply it to Question \ref{T:conj2}.

We start by showing how a normal conditional expectation $E:\mathcal{M} \to \mathcal{N}$ induces an inclusion of $L^p$ spaces.  Since $L^p(s(E)\mathcal{M}s(E)) \subset L^p(\mathcal{M})$ naturally, we may assume that $E$ is faithful.  One method is by Kosaki's adaptation of the complex interpolation method [K1].  Assume that $\mathcal{N}$ is $\sigma$-finite, fix a faithful state $\varphi \in \mathcal{N}_*$, and consider the left embedding of $\mathcal{N}$ in $\mathcal{N}_*$: $x \mapsto x\varphi$.  Then $L^p(\mathcal{N})$ arises as the interpolated Banach space at $1/p$ ([K1], Theorem 9.1); more precisely, we have
\begin{equation} \label{E:interp}
L^p(\mathcal{N}) \varphi^{1/q} = [\mathcal{N}, \mathcal{N}_*]_{1/p} \simeq L^p(\mathcal{N}), \qquad 1/p + 1/q = 1.
\end{equation}
Here the equality is meant as \textit{sets}, while the isomorphism is an isometric identification of Banach spaces.

We isometrically include this interpolation couple, by $E_*$, in the interpolation couple for $(\mathcal{M}, \mathcal{M}_*)$ arising from the embedding $y \mapsto y(\varphi \circ E)$.  The reader can check that the following diagram commutes, with the horizontal compositions being identity maps:
\[
\begin{CD}
\mathcal{N} @>>> \mathcal{M} @>E>> \mathcal{N}\\
@VVV            @VVV                  @VVV\\
\mathcal{N}_*  @>E_*>>  \mathcal{M}_* @>\text{restriction}>> \mathcal{N}_*
\end{CD}
\]
It is important here that $E$ is a bimodule map: $E(n_1 m n_2) = n_1 E(m) n_2$!  Then by general interpolation theory (e.g. Theorem 1.2 in [K1]) one may interpolate these 1-complemented inclusions to get a 1-complemented inclusion at the $L^p$ level.  By \eqref{E:interp} we know that the map is densely defined by $x \varphi^{1/p} \mapsto x (\varphi \circ E)^{1/p}$.

\begin{proposition} \label{T:indt}
The $L^p$ isometry constructed in the previous paragraph is independent of the choice of $\varphi$.
\end{proposition}

\begin{proof}
We show that when there is some $C < \infty$ so that $\varphi^{2/p} \le C \psi^{2/p}$ as operators in the modular algebra,
\begin{equation} \label{E:opwt}
x \varphi^{1/p} = y \psi^{1/p} \Rightarrow x(\varphi \circ E)^{1/p} = y (\psi \circ E)^{1/p}, \qquad x,y \in \mathcal{N}.
\end{equation}
This is sufficient, because (see Section 1 of [JS])
$$ \varphi^{2/p} \le C \psi^{2/p} \Rightarrow \varphi^{1/p} = z \psi^{1/p} \text{ for some }z \in \mathcal{N},$$
and \eqref{E:opwt} then implies that the embeddings for $\varphi$ and $\psi$ agree on the dense set $\{x\varphi^{1/p} \mid x \in \mathcal{N}\}.$  For any two faithful normal states $\varphi, \theta$, we can conclude that the embeddings each equal the embedding for $\frac{(\varphi^{2/p} + \theta^{2/p})^{p/2} }{\| (\varphi^{2/p} + \theta^{2/p})^{p/2} \|}$ and therefore equal each other.

To prove \eqref{E:opwt}, first recall the Connes cocycle derivative equation (Corollary VIII.4.22 of [T2])
\begin{equation} \label{E:cocycle}
(D(\varphi \circ E) : D(\psi \circ E))_t = (D\varphi: D\psi)_t.
\end{equation}
The condition $\varphi^{2/p} \le C \psi^{2/p}$ guarantees that $(D\varphi: D\psi)_t$ extends to a continuous $\mathcal{N}$-valued function on the strip $\{0 \ge \text{Im } z \ge -1/p \}$ which is analytic on the interior ( $\sim$ Theorem VIII.3.17 of [T2], or see [S2]).  We calculate
\begin{align}
\notag x \varphi^{1/p} = y \psi^{1/p} &\Rightarrow x (D\varphi: D\psi)_{-i/p} = y \\
\notag &\Rightarrow x(D(\varphi \circ E) : D(\psi \circ E))_{-i/p} = y \\
\notag &\Rightarrow x(\varphi \circ E)^{1/p} = y (\psi \circ E)^{1/p}.
\end{align}

\end{proof}

So we can avoid interpolation (and the choice of a reference state) altogether by simply writing the $L^p$ isometry as
\begin{equation} \label{E:canon}
\varphi^{1/p} \mapsto (\varphi \circ E)^{1/p}, \qquad \varphi \in \mathcal{N}_*^+,
\end{equation}
and extending linearly off the positive cone.  In a manner similar to the proof of Proposition \ref{T:indt}, one can show that the family of inclusions \eqref{E:canon} preserves the duality between $L^p(\mathcal{N})$ and $L^q(\mathcal{N})$ ($1/p + 1/q = 1$).  It is also worth noting that right-hand embeddings of the form $\mathcal{N} \ni x \mapsto \varphi x \in \mathcal{N}_*$, or even others, will necessarily produce the same $L^p$ isometry, namely \eqref{E:canon}.

In case $\mathcal{N}$ is not $\sigma$-finite, \eqref{E:canon} defines an $L^p$ isometry on each $q \mathcal{N} q$, where $q$ is a $\sigma$-finite projection in $\mathcal{N}$.  Every finite set of vectors in $L^p(\mathcal{N})$ belongs to some such $qL^p(\mathcal{N})q$, as the left and right supports of each vector belong to the lattice of $\sigma$-finite projections.  Being defined by \eqref{E:canon}, these $L^p$ isometries agree on their common domains and so define a global $L^p$ isometry.

Suppressed in the above scenario is an inclusion map $\iota: \mathcal{N} \hookrightarrow \mathcal{M}$, which is of course multiplicative.  What if it is antimultiplicative?  Since we have already discussed the effect of the condtional expectation, let us consider only the map on $L^p(\mathcal{M})$ induced by a normal *-antiautomorphism $\alpha: \mathcal{M} \to \mathcal{M}$.

Fixing faithful $\varphi$, we set up the interpolation as follows.  In the domain, $\mathcal{M} \subset \mathcal{M}_*$ via $x \hookrightarrow x\varphi$, while in the range, $\mathcal{M} \subset \mathcal{M}_*$ via $y \hookrightarrow (\varphi \circ \alpha^{-1})y$.  This gives us the commutative diagram
\[
\begin{CD}
\mathcal{M} @>\alpha>> \mathcal{M}\\
@VVV            @VVV\\
\mathcal{M}_*  @>>(\alpha^{-1})_*>  \mathcal{N}_*
\end{CD}
\]
Because the inclusions are commuting and surjective, we get a surjective $L^p$ isometry which by \eqref{E:interp} is densely defined by $x\varphi^{1/p} \mapsto (\varphi \circ \alpha^{-1})^{1/p} \alpha(x)$.

Once again this map is independent of the choice of $\varphi$.  The proof is the same as that of Proposition \ref{T:indt}: it suffices to verify
\begin{equation} \label{E:flip}
x \varphi^{1/p} = y \psi^{1/p} \Rightarrow (\varphi \circ \alpha^{-1})^{1/p} \alpha(x) = (\psi \circ \alpha^{-1})^{1/p} \alpha(y)
\end{equation}
under the assumption $\varphi^{2/p} \le \psi^{2/p}$.
Temporarily assuming the cocycle identity
\begin{equation} \label{E:anticocyc}
(D(\psi \circ \alpha^{-1}): D(\varphi \circ \alpha^{-1}))_t = \alpha((D\varphi:D\psi)_{-t}),
\end{equation}
we have
\begin{align}
\notag x \varphi^{1/p} = y \psi^{1/p} &\Rightarrow x (D\varphi: D\psi)_{-i/p} = y \\
\notag &\Rightarrow \alpha[(D\varphi : D\psi)_{-i/p}] \alpha(x) = \alpha(y) \\
\notag &\Rightarrow [(D(\psi \circ \alpha^{-1}): D(\varphi \circ \alpha^{-1}))_{i/p}] \alpha(x) = \alpha(y) \\
\notag &\Rightarrow (\varphi \circ \alpha^{-1})^{1/p} \alpha(x) = (\psi \circ \alpha^{-1})^{1/p} \alpha(y).
\end{align}
Of course it remains to show

\begin{lemma}
Equation \eqref{E:anticocyc} holds.
\end{lemma}

\begin{proof}
Let $\alpha: \mathcal{M} \to \mathcal{M}$ be a normal *-antiautomorphism.  We first claim that
\begin{equation} \label{E:modaut}
\sigma_t^{\varphi \circ \alpha^{-1}} = \alpha \circ \sigma_{-t}^{\varphi} \circ \alpha^{-1}.
\end{equation}
Recall ([T2], Theorem VIII.1.2) that the modular automorphism group for $\varphi \circ \alpha^{-1}$ is the unique one-parameter automorphism group which (1) is $\varphi \circ \alpha^{-1}$ invariant and (2) satisfies the KMS condition for $\varphi \circ \alpha^{-1}$.  We check that the right-hand side above meets the conditions.  For the first,
$$\varphi \circ \alpha^{-1}(\alpha \circ \sigma_{-t}^{\varphi} \circ \alpha^{-1} (x)) = \varphi \circ \sigma_{-t}^{\varphi} \circ \alpha^{-1} (x) = \varphi \circ \alpha^{-1}(x), \qquad x \in \mathcal{M}.$$
For the second, fix $x,y \in \mathcal{M}$.  Use the KMS condition for $\varphi$ to find a function $F = F_{\alpha^{-1}(x), \alpha^{-1}(y)}$ on $\{ 0 \le \text{Im }z \le 1\}$ which is bounded, continuous, and analytic on the interior - these properties are assumed but not stated for later functions on the strip - with boundary values
$$F(t) = \varphi[(\sigma_t^\varphi(\alpha^{-1}(x))) (\alpha^{-1}(y))], \qquad F(t+i) = \varphi[(\alpha^{-1}(y)) (\sigma_t^\varphi (\alpha^{-1}(x)))].$$
Notice
$$F(t) = \varphi \circ \alpha^{-1}[(y)(\alpha \circ \sigma_t^\varphi \circ \alpha^{-1}(x))], \qquad F(t+i) = \varphi \circ \alpha^{-1}[(\alpha \circ \sigma_t^\varphi \circ \alpha^{-1}(x)) (y)].$$
Then $G(z) \triangleq F(i-z)$ is a function on the strip which satisfies 
$$G(t) = \varphi \circ \alpha^{-1}[(\alpha \circ \sigma_{-t}^\varphi \circ \alpha^{-1}(x))( y)], \qquad G(t+i) = \varphi \circ \alpha^{-1}[(y) (\alpha \circ \sigma_{-t}^\varphi \circ \alpha^{-1}(x))].$$
Thus $\alpha \circ \sigma_{-t}^{\varphi} \circ \alpha^{-1}$ satisfies the KMS condition for $\varphi \circ \alpha^{-1}$ at any $x,y$, and \eqref{E:modaut} is proved.

We establish \eqref{E:anticocyc} in a similar way.  Choose $x,y \in \mathcal{M}$, and by ([T2], Theorem VIII.3.3) find a function $F = F_{\alpha^{-1}(x), \alpha^{-1}(y)}$ on the strip with
$$F(t) = \varphi[(D\varphi:D\psi)_t \sigma_t^\psi( \alpha^{-1}(y)) \alpha^{-1}(x)],$$
$$F(t+i) = \psi[\alpha^{-1}(x)(D\varphi:D\psi)_t \sigma_t^\psi(\alpha^{-1}(y))].$$
We rewrite this, using cocycle relations and \eqref{E:modaut}:
\begin{align}
\notag F(t) &= \varphi[\sigma_t^\varphi(\alpha^{-1}(y)) (D\varphi:D\psi)_t \alpha^{-1}(x)] \\
\notag &= \varphi \circ \alpha^{-1}[ x \: \alpha((D\varphi:D\psi)_t) (\alpha \circ \sigma_t^\varphi \circ \alpha^{-1}(y))] \\
\notag &= \varphi \circ \alpha^{-1}[ x \: \alpha((D\varphi:D\psi)_t) \sigma_{-t}^{\varphi \circ \alpha^{-1}}(y)].
\end{align}
Analogously, we obtain
$$F(t + i) = \psi \circ \alpha^{-1}[\alpha((D\varphi:D\psi)_t) (\sigma_{-t}^{\varphi \circ \alpha^{-1}}(y)) x].$$
Then $G(z) \triangleq F(i - z)$ is a function on the strip which satisfies
$$G(t) = \psi \circ \alpha^{-1}[\alpha((D\varphi:D\psi)_{-t}) (\sigma_t^{\varphi \circ \alpha^{-1}}(y)) x],$$ $$ G(t+i) = \varphi \circ \alpha^{-1}[ x \: \alpha((D\varphi:D\psi)_{-t}) \sigma_t^{\varphi \circ \alpha^{-1}}(y)].$$
Again by Theorem VIII.3.3 of [T2], the existence of such a $G$ for any $x,y$ implies that $\alpha((D\varphi:D\psi)_{-t}) = (D(\psi \circ \alpha^{-1}): D(\varphi \circ \alpha^{-1}))_t$, so we are done.
\end{proof}

The non-$\sigma$-finite case can be handled as in our previous discussion.  We summarize these observations in
\begin{proposition} \label{T:lpcond}
Let $\alpha$ be a normal *-isomorphism or *-antiisomorphism from $\mathcal{M}_1$ into $\mathcal{M}_2$, and suppose that there is a normal conditional expectation $E:\mathcal{M}_2 \to \alpha(\mathcal{M}_1).$  Then the map
$$\varphi^{1/p} \mapsto (\varphi \circ \alpha^{-1} \circ E)^{1/p}, \qquad \varphi \in (\mathcal{M}_1)_*^+,$$
extends off the positive cone to a (typical) isometry $L^p(\mathcal{M}_1) \to L^p(\mathcal{M}_2)$.
\end{proposition}

\section{Modular theory and projections onto Jordan subalgebras} \label{S:jordan}

What happens to Proposition \ref{T:lpcond} when $\alpha$ and $E$ are replaced with a normal Jordan *-monomorphism $J$ and a normal positive projection $P$?  To answer this, we first recall the Jordan version of modular theory.  We can still construct $L^p$ isometries by interpolation, but the lack of a \textit{relative} modular theory for Jordan algebras (along the lines of \eqref{E:cocycle} and \eqref{E:anticocyc}) has prevented us from concluding in general that there is no dependence on the choice of reference state.  In Section \ref{S:factor} we will introduce hypotheses which remove this (possible) dependence.

\bigskip

We will need to compare some modular objects for linear functionals on $\mathcal{M}_1$, $J(\mathcal{M}_1)$, and $\mathcal{M}_2$.  Since the second is usually not a von Neumann algebra, we must use in place of a modular automorphism group the \textit{modular cosine family} [HH].  This is defined generally for a normal (faithful) state $\varphi$ on a JBW-algebra $\mathcal{N}$; we do not need the generality but do need the five conditions which uniquely characterize $\rho^\varphi_t$ ([HH], Theorem 3.3):
\begin{enumerate}
\item each $\rho^\varphi_t$ is a positive normal unital linear map;
\item $\mathbb{R} \ni t \mapsto \rho^\varphi_t(x)$ is $\sigma$-weakly continuous for each $x \in \mathcal{N}$;
\item $\rho^\varphi_0 = \text{id}_\mathcal{N}$ and $\rho^\varphi_s \circ \rho^\varphi_t = (1/2)[\rho^\varphi_{s+t} + \rho^\varphi_{s-t}]$;
\item $\varphi(\rho^\varphi_t(a) \bullet b) = \varphi(a \bullet \rho^\varphi_t(b))$;
\item the sesquilinear form
$$s_\varphi(a,b) = \int_{-\infty}^\infty \varphi(\rho^\varphi_t(a) \bullet b^*) (\cosh \pi t)^{-1} dt$$
is self-polar.
\end{enumerate}
(A positive sesquilinear form $s(\cdot, \cdot)$ is said to be \textit{self-polar} ([HH], [Wo]) if (i) $s(a,b) \ge 0, \: \forall a,b \in \mathcal{N}_+$; and (ii) the set of linear functionals $\{s(\cdot, h) \mid 0 \le h \le 1\}$ is weak*-dense in $\{\psi \in \mathcal{N}^*_+ \mid \psi \le \varphi \}$.)  When $\mathcal{N}$ is a von Neumann algebra, we have
$$\rho^\varphi_t = (1/2)(\sigma^\varphi_t + \sigma^\varphi_{-t})$$
and
\begin{equation} \label{E:selfpolar}
s_\varphi(a,b) = <\varphi^{1/4}a \varphi^{1/4},\: \varphi^{1/4} b \varphi^{1/4}>.
\end{equation}

The following result of Haagerup and St\o rmer is the Jordan algebra version of Takesaki's theorem [T1] for conditional expectations.  We specialize it to our situation.
\begin{theorem} \label{T:hs} ([HS2], Theorem 4.2)
Let $J: \mathcal{M}_1 \to \mathcal{M}_2$ be a normal Jordan *-monomorphism.  Let $\psi \in (\mathcal{M}_2)_*$ be a faithful state, and denote by $\theta$ the restriction $\psi \mid_{J(\mathcal{M}_1)}$.  Then the following are equivalent:
\begin{enumerate}
\item There is a normal positive projection $P: \mathcal{M}_2 \to J(\mathcal{M}_1)$ such that $\psi = \theta \circ P$;
\item $s_\theta = s_\psi \mid_{J(\mathcal{M}_1) \times J(\mathcal{M}_1)}$;
\item $\rho_t^\theta = \rho_t^\psi \mid_{J(\mathcal{M}_1)}, \: \forall t \in \mathbb{R}$.
\end{enumerate}
\end{theorem}
When these conditions hold, $P$ can be defined by
\begin{equation} \label{E:spolar}
s_\psi(y,J(x)) = s_\psi(P(y), J(x)), \qquad x \in \mathcal{M}_1, \: y \in \mathcal{M}_2.
\end{equation}
Note that the analogue for a normal $\psi$-preserving condtional expectation $E:\mathcal{N} \to \mathcal{M}$ is
\begin{equation} \label{E:cexp}
\psi(yx) = \psi(E(y)x), \qquad y \in \mathcal{N}, \: x \in \mathcal{M}.
\end{equation}

\bigskip

For a faithful normal state $\varphi$, we will make use of the following transform, familiar from Tomita-Takesaki theory:
\begin{equation} \label{E:Phi1}
\Phi_\varphi(x) = \int_{-\infty}^\infty \sigma_t^\varphi(x) (\cosh \pi t)^{-1} dt = \int_{-\infty}^\infty \rho_t^\varphi(x) (\cosh \pi t)^{-1} dt.
\end{equation}
We have the equality ([vD], Lemma 4.1)
\begin{equation} \label{E:Phi2}
\varphi^{1/2} x \varphi^{1/2} = (1/2)[ \Phi_\varphi(x) \varphi + \varphi \Phi_\varphi(x)] \triangleq \varphi \bullet \Phi_\varphi(x).
\end{equation}
So $\Phi_\varphi(x)$ is the Jordan derivative of $\varphi^{1/2} x \varphi^{1/2}$ with respect to $\varphi$.  (In fact $\Phi_\varphi$ is a right inverse for $\rho_{i/2}^\varphi$ (which makes sense for a dense set, via analytic continuation).  One might also observe that \eqref{E:Phi1} and \eqref{E:Phi2} prove the implication (3) $\to$ (2) of Theorem \ref{T:hs}, as $\varphi^{1/2} x \varphi^{1/2}(y) = s_\varphi(x,y^*)$.)

\bigskip

We will assume that $\mathcal{M}_1$ is $\sigma$-finite, and for the time being, assume also that $P$ is faithful.  Choose a faithful state $\varphi \in (\mathcal{M}_1)_*$.  Now we set up Kosaki's complex interpolation method [K1] again, but we require the symmetric embedding
\begin{equation} \label{E:central}
\iota_1: \mathcal{M}_1 \hookrightarrow (\mathcal{M}_1)_*, \qquad x \mapsto \varphi^{1/2} x \varphi^{1/2}. 
\end{equation}
This gives us the interpolation spaces 
\begin{equation} \label{E:kosaki}
\varphi^{1/2q} L^p(\mathcal{M}_1) \varphi^{1/2q} = [\mathcal{M}_1, (\mathcal{M}_1)_*]_{1/p} \simeq L^p(\mathcal{M}_1), \qquad 1/p + 1/q = 1,
\end{equation}
where again the equality is between sets and the isomorphism is isometric.
Accordingly, we can base a construction of $L^p(\mathcal{M}_2)$ on a symmetric embedding $\iota_2$ using the faithful state $\varphi \circ J^{-1} \circ P$.  We will often compress notation by denoting the image of $(J^{-1} \circ P)_*$ by a bar, so $(\varphi \circ J^{-1} \circ P) = \bar{\varphi}$.

\begin{theorem} \label{T:compatible}
With the setup of the preceding paragraph, the following diagram commutes:
\[
\begin{CD}
\mathcal{M}_1 @>J>> \mathcal{M}_2 @>J^{-1} \circ P>> \mathcal{M}_1\\
@V\iota_1VV            @V\iota_2VV                  @V\iota_1VV\\
(\mathcal{M}_1)_*  @>(J^{-1} \circ P)_*>> (\mathcal{M}_2)_* @>J_*>> (\mathcal{M}_1)_*
\end{CD}
\]
This induces the 1-complemented inclusion of $L^p$ spaces
$$ L^p(\mathcal{M}_1) \hookrightarrow L^p(\mathcal{M}_2) \twoheadrightarrow L^p(\mathcal{M}_1),$$
densely defined by
\begin{equation} \label{E:ident}
\varphi^{1/2p}x \varphi^{1/2p} \to \bar{\varphi}^{1/2p} J(x) \bar{\varphi}^{1/2p}, \qquad x \in \mathcal{M}_1,
\end{equation}
and
\begin{equation} \label{E:lpproj}
\bar{\varphi}^{1/2p} y \bar{\varphi}^{1/2p} \to \varphi^{1/2p} (J^{-1} \circ P(y)) \varphi^{1/2p}, \qquad y \in \mathcal{M}_2.
\end{equation}
\end{theorem}

\begin{proof}
Theorem \ref{T:hs} tells us that the restriction of $\rho_t^{(\varphi \circ J^{-1} \circ P)}$ to $J(\mathcal{M}_1)$ is $\rho_t^{(\varphi \circ J^{-1}) }$.  By consulting the five conditions characterizing $\rho_t^{(\varphi \circ J^{-1})}$, one checks that it agrees with $J \circ \rho_t^\varphi \circ J^{-1}$.  (This ``must" be true, as $J$ is a normal Jordan *-isomorphism onto its image.)  These facts imply
\begin{align*}
\Phi_{\bar{\varphi}}(J(x)) &= \int \rho_t^{ \bar{\varphi}}(J(x)) (\cosh \pi t)^{-1} dt \\
&= \int \rho_t^{(\varphi \circ J^{-1})}(J(x)) (\cosh \pi t)^{-1} dt \\
&= \int J \circ \rho_t^\varphi \circ J^{-1}(J(x)) (\cosh \pi t)^{-1} dt\\
&= J\left[\int \rho_t^\varphi (x) (\cosh \pi t)^{-1} dt\right] \\
&= J(\Phi_\varphi(x)).
\end{align*}

Now we are ready to check that the diagram in Theorem \ref{T:compatible} commutes, starting with the left-hand square.  This amounts to showing that for $x \in \mathcal{M}_1, \: y \in \mathcal{M}_2$,
\begin{equation} \label{E:compat}
\bar{\varphi}^{1/2} J(x) \bar{\varphi}^{1/2} (y) = [(\varphi^{1/2} x \varphi^{1/2}) \circ J^{-1} \circ P] (y).
\end{equation}
We calculate
\begin{align*}
\bar{\varphi}^{1/2} J(x) \bar{\varphi}^{1/2} (y) &= [\bar{\varphi} \bullet \Phi_{\bar{\varphi}} (J(x))] (y) \\
&= \bar{\varphi} (\Phi_{\bar{\varphi}}(J(x)) \bullet y) \\
&= \bar{\varphi} ((J(\Phi_{\varphi}(x))) \bullet y) \\
&\overset{*}{=} \varphi (\Phi_\varphi(x) \bullet (J^{-1} \circ P(y))) \\
&= \varphi^{1/2} x \varphi^{1/2}(J^{-1} \circ P(y)) \\
&= [(\varphi^{1/2} x \varphi^{1/2}) \circ J^{-1} \circ P](y).
\end{align*}
We used Lemma \ref{T:stormer}(2) in the equality marked $\overset{*}{=}$.

For the right-hand square, we need to demonstrate that for $x \in \mathcal{M}_1,\: y \in \mathcal{M}_2$,
$$[(\bar{\varphi}^{1/2} y \bar{\varphi}^{1/2}) \circ J](x) = (\varphi^{1/2} (J^{-1} \circ P(y)) \varphi^{1/2})(x).$$
But this equation is equivalent to \eqref{E:compat}, which was just shown.

It again follows from general interpolation theory that the inclusion and norm one projection extend to the interpolated spaces.  Since $\varphi^{1/2} x \varphi^{1/2}$ is identified with $ \bar{\varphi}^{1/2}J(x) \bar{\varphi}^{1/2}$, the equality \eqref{E:kosaki} gives us \eqref{E:ident} and \eqref{E:lpproj}.

\end{proof}

\textit{Remark 1.}  Theorem \ref{T:compatible} holds without change if $P$ is not faithful.  With $q$ the support of $P$, replace $\mathcal{M}_2$ by $q\mathcal{M}_2q$, and notice that $L^p(q\mathcal{M}_2 q) \simeq qL^p(\mathcal{M}_2)q \subset L^p(\mathcal{M}_2)$.

It also seems possible, if not pleasant, to extend Theorem \ref{T:compatible} to non-$\sigma$-finite $\mathcal{M}_1$ by using a weight.  When $P$ is a normal conditional expectation, most of the necessary tools are in [Te2] and [I].

\textit{Remark 2.}  As mentioned, typicality may be related to a Jordan version of \eqref{E:cocycle} and \eqref{E:anticocyc}.  Cocycles are not symmetric objects (the ``handedness" is apparent in the modular algebra realization $(D\varphi: D\psi)_t = \varphi^{it} \psi^{-it}$), so one cannot expect either of these equations for the projection $P$.  We can derive at least something analogous using \eqref{E:compat}.  Let faithful $\varphi \le C \psi \in (\mathcal{M}_1)^+_*$ for some $C < \infty$, so $y = (D\varphi: D\psi)_{-i/2} (=\varphi^{1/2} \psi^{-1/2})$ exists in $\mathcal{M}_1$.  Now write
$$\bar{\varphi} = \varphi \circ J^{-1} \circ P = \psi^{1/2} y^* y \psi^{1/2} \circ J^{-1} \circ P =\bar{\psi}^{1/2} J(y^* y) \bar{\psi}^{1/2}.$$
If we write $z = (D\bar{\varphi}: D\bar{\psi})_{-i/2}$, this gives $z^*z = J(y^* y)$.  Taking square roots, $|z| = J(|y|)$, or
$$|(D\bar{\varphi}: D\bar{\psi})_{-i/2}| = J(|(D\varphi: D\psi)_{-i/2}|).$$

\section{Factorization and typical isometries} \label{S:factor}

We have not been able to show that the isometry constructed in \eqref{E:ident} is generally independent of the choice of $\varphi$.  We now introduce a hypothesis which removes the (possible) dependence.

\begin{proposition} \label{T:factor}
Suppose $J: \mathcal{M}_1 \to \mathcal{M}_2$ is a normal Jordan *-monomorphism and $P: \mathcal{M}_2 \to J(\mathcal{M}_1)$ is a normal positive projection, with $P$ factoring as
\begin{equation} \label{E:factor}
P = P' \circ F, \qquad F: \mathcal{M}_2 \to J(\mathcal{M}_1)'' \text{ a normal conditional expectation.}
\end{equation}
Then \eqref{E:typical} (taking $w=1$) defines an $L^p$ isometry.  If $\mathcal{M}_1$ is $\sigma$-finite, then this agrees with the $L^p$ isometry \eqref{E:ident} obtained in Theorem \ref{T:compatible} for any faithful normal state $\varphi \in (\mathcal{M}_1)_*$.
\end{proposition}

\begin{proof}
It suffices to prove these statements for the map
\begin{equation} \label{E:double}
L^p(\mathcal{M}_1)_+ \ni \varphi^{1/p} \mapsto (\varphi \circ J^{-1} \circ P')^{1/p} \in L^p(J(\mathcal{M}_1)'')_+,
\end{equation}
since precomposition with the conditional expectation $F$ embeds $L^p(J(\mathcal{M}_1)'')$ in $L^p(\mathcal{M}_2)$ as discussed in Section \ref{S:conditioned}.

Now $J$ is the direct sum of a *-homomorphism $\pi$ and a *-antiisomorphism $\pi'$, so $\pi(1)$ and $\pi'(1)$ are orthogonal and central in $J(\mathcal{M}_1)''$.  We may assume that the abelian summand of $\mathcal{M}_1$, if it exists, is in the support of $\pi$ and not $\pi'$.  This decomposes $\mathcal{M}_1$ into three summands: $s(\pi)(1-s(\pi')) \mathcal{M}_1$, $(1-s(\pi))s(\pi') \mathcal{M}_1$, and $s(\pi)s(\pi') \mathcal{M}_1$.  Note that the images of the first two are central summands of $J(\mathcal{M}_1)$ which are multiplicatively closed, so they are also central summands of $J(\mathcal{M}_1)''$.  Since $P$ and $P'$ respect this central decomposition (by Lemma \ref{T:stormer}), $P'$ restricts to the identity on these, and we are in the situation of Proposition \ref{T:lpcond}.  It is left to discuss the third summand, so in the remainder of the proof we assume that $s(\pi) = s(\pi') = 1$ and that $\mathcal{M}_1$ has no abelian summand.

We first claim
$$J(\mathcal{M}_1)'' = \{\pi(x) \oplus \pi'(y) \mid x,y \in \mathcal{M}_1 \} \simeq \mathcal{M}_1 \oplus \mathcal{M}_1^{\text{op}}.$$
Let $v \in \mathcal{M}_1$ be a partial isometry with $vv^* = p \perp q = v^*v$.  Then $J(\mathcal{M}_1)''$ contains
$$(\pi(v) \oplus \pi'(v))(\pi(p) \oplus \pi'(p)) = 0 \oplus \pi'(v)$$
and is closed under multiplication, addition, and adjoints, whence it is easy to verify the claim.

By Lemma \ref{T:stormer}(4), $P'(\pi(1) \oplus 0) = J(\lambda) \in J(\mathcal{Z}(\mathcal{M}_1))$, where $0 \le \lambda \le 1$.  Then
\begin{align} \label{E:lambda}
P'(\pi(x) \oplus 0) &= P'[(\pi(x) \oplus \pi'(x)) \bullet (\pi(1) \oplus 0)] \\
\notag &= (\pi(x) \oplus \pi'(x)) \bullet P'(\pi(1) \oplus 0) \\
\notag &= \pi(x)\pi(\lambda) \oplus \pi'(x)\pi'(\lambda)\\
\notag &= J(x\lambda),
\end{align}
and similarly
$$P'(0 \oplus \pi'(y)) = \pi(y)\pi(1 - \lambda) \oplus \pi'(y)\pi'(1 - \lambda) = J(y(1 - \lambda)).$$
Note that $\lambda$ and $1-\lambda$ must be nonsingular, because faithfulness of $P'$ is a consequence of $s(P) = P(1) (= J(1))$.  We obtain that
$$\varphi \circ J^{-1} \circ P' = (\lambda \varphi \circ \pi^{-1}) \oplus ((1-\lambda) \varphi \circ \pi'^{-1}).$$
The map \eqref{E:double} is then 
$$\varphi^{1/p} \mapsto (\lambda \varphi \circ \pi^{-1})^{1/p} \oplus ((1-\lambda)\varphi \circ \pi'^{-1})^{1/p}.$$
Since $\pi$ and $\pi'$ induce (surjective) isometric isomorphisms at the $L^p$ level in each summand, this map has a linear extension to all of $L^p(\mathcal{M}_1)$.  The image of $\xi \in L^p(\mathcal{M}_1)$ is a vector whose two orthogonal components have norms $\|\lambda^{1/p}  \xi \|$ and $\|(1- \lambda)^{1/p} \xi \|$; it has total norm
$$(\|\lambda |\xi|^p \| + \|(1 - \lambda) |\xi|^p \|)^{1/p} = \| |\xi|^p \|^{1/p} = \| \xi \|.$$
Therefore the map is isometric.

The second assertion of the proposition almost follows from the discussion in Section \ref{S:conditioned} of $L^p$ isometries induced by (possibly antimultiplicative) 1-complemented inclusions.  There we noted that any left or right embedding gave the same isometry, and so was typical.  Here we make this statement explicit for symmetric embeddings.  So suppose $\varphi^{1/2p}x\varphi^{1/2p} = \psi^{1/2p} y \psi^{1/2p} \in L^p(\mathcal{M}_1)_+$.  By taking square roots, we can find a partial isometry $v \in \mathcal{M}_1$ with $vx^{1/2}\varphi^{1/2p} = y^{1/2} \psi^{1/2p}$.  Considering first the antimultiplicative embedding, we know by \eqref{E:flip} that
$$(\varphi \circ \pi'^{-1})^{1/2p}\pi'(vx^{1/2}) =  (\psi\circ \pi'^{-1})^{1/2p}\pi'(y^{1/2}),$$
so
\begin{align*}
(\varphi \circ \pi'^{-1})^{1/2p} \pi'(x) &(\varphi \circ \pi'^{-1})^{1/2p}\\
&= [(\varphi \circ \pi'^{-1})^{1/2p}\pi'(vx^{1/2})] [(\varphi \circ \pi'^{-1})^{1/2p}\pi'(vx^{1/2})]^* \\
&= [(\psi\circ \pi'^{-1})^{1/2p}\pi'(y^{1/2})] [(\psi\circ \pi'^{-1})^{1/2p}\pi'(y^{1/2})]^* \\
&= (\psi\circ \pi'^{-1})^{1/2p}\pi'(y) (\psi\circ \pi'^{-1})^{1/2p}.
\end{align*}
Obviously this relation extends off the positive cone.  A similar calculation holds for $\pi$, so it holds for $J$.  This means that the $L^p$ isometry \eqref{E:ident} does not depend on the choice of state and therefore agrees with \eqref{E:double}.
\end{proof}

The papers [HS2], [St3] consider exactly the factorization \eqref{E:factor} for projections onto arbitrary JW-subalgebras, which includes our situation.  They do conclude the factorization in case
\begin{enumerate}
\item $\mathcal{Z}(J(\mathcal{M}_1)) = \mathcal{Z}(J(\mathcal{M}_1)'')$; or 
\item $(\mathcal{M}_1, \tau)$ is semifinite and the weight $\bar{\tau} = \tau \circ J^{-1} \circ P$ is semifinite.
\end{enumerate}
In our situation, condition (1) precludes the presence of both multiplicative and antimultiplicative homomorphisms on a non-abelian summand, so $J(\mathcal{M}_1) = J(\mathcal{M}_1)''$ and the conclusion is automatic.  Our final result subsumes condition (2) by assuming only that $\mathcal{M}_1$ is AS.

\begin{proposition} \label{T:asfactor}
Let $J:\mathcal{M}_1 \to \mathcal{M}_2$ be a normal Jordan *-monomorphism, and $P:\mathcal{M}_2 \to J(\mathcal{M}_1)$ be a normal positive projection.  If $\mathcal{M}_1$ is AS, then $P$ factors as $P \mid_{J(\mathcal{M}_1)''} \circ F$, where $F: \mathcal{M}_2 \to J(\mathcal{M}_1)''$ is a normal conditional expectation.
\end{proposition}

\begin{proof} First note that by the arguments in the beginning of the proof of Proposition \ref{T:factor}, we may assume that $J = \pi \oplus \pi'$, where $\mathcal{M}_1$ has no abelian summand and both $\pi$, $\pi'$ are faithful.

We start with the case where $(\mathcal{M}_1, \tau)$ is finite, writing $\bar{\tau} = \tau \circ J^{-1} \circ P$.  In the proof of Theorem \ref{T:cgyt} we saw that $J(\mathcal{M}_1)''$ is pointwise invariant under $\sigma^{\bar{\tau}}$, so by Takesaki's theorem [T1] there is a normal $\bar{\tau}$-preserving conditional expectation $F: \mathcal{M}_2 \to J(\mathcal{M}_1)''$.  Now we use equations \eqref{E:spolar} and \eqref{E:cexp} to calculate
$$s_{\bar{\tau}}(P(F(y)), J(x)) = s_{\bar{\tau}}(F(y), J(x)) = \bar{\tau}^{1/2} J(x^*) \bar{\tau}^{1/2}(F(y)) = \bar{\tau}(J(x^*)F(y))$$ $$ = \bar{\tau}(J(x^*)y) = \bar{\tau}^{1/2} J(x^*) \bar{\tau}^{1/2}(y) = s_{\bar{\tau}}(y, J(x)) = s_{\bar{\tau}}(P(y), J(x))$$
for any $x \in \mathcal{M}_1,\: y \in \mathcal{M}_2$.  This implies that $P \circ F = P$, and the finite case is complete.

Now let $\{\mathcal{M}_\alpha, E_\alpha\}$ be a paving of $\mathcal{M}_1$ by finite subalgebras, and denote by $1_\alpha$ the unit of $\mathcal{M}_\alpha$ (which may not be the unit of $\mathcal{M}_1$).  Write 
$$\bar{E}_\alpha = J \circ E_\alpha \circ J^{-1}: J(\mathcal{M}_1) \to J(\mathcal{M}_\alpha).$$
Then $\bar{E}_\alpha \circ P$ is a positive normal projection onto the Jordan image of a finite algebra, so by the first part of the proof we can factor it as
\[
\begin{CD}
\mathcal{M}_2 @>P>> J(\mathcal{M}_1)\\
@VF_{\alpha}VV            @VV\bar{E}_{\alpha}V \\
J(\mathcal{M}_\alpha)''  @>S_{\lambda_\alpha}>> J(\mathcal{M}_\alpha)
\end{CD}
\]
Here $F_\alpha$ is a normal conditional expectation, and $S_{\lambda_\alpha}$ is the symmetrization guaranteed by \eqref{E:lambda},
\begin{equation} \label{E:lamb}
S_{\lambda_\alpha}: \pi(x) \oplus \pi'(y) \mapsto J(\lambda_\alpha x + (1-\lambda_\alpha) y),
\end{equation}
where $0 \le \lambda_\alpha \le 1_\alpha$ is an element of $J(\mathcal{Z}(\mathcal{M}_\alpha))$.  We also have, as before, that the nonsingular element $J(\lambda) \triangleq P(\pi(1))$ is in $J(\mathcal{Z}(\mathcal{M}_1)).$  The commuting square gives us the relation
$$\bar{E}_\alpha (J(\lambda)) = \bar{E}_\alpha \circ P(\pi(1)) = S_{\lambda_\alpha} \circ F_\alpha (\pi(1)) = S_{\lambda_\alpha} (\pi(1_\alpha)) = J(\lambda_\alpha).$$
This implies that $\lambda_\alpha = E_\alpha(\lambda) \overset{s}{\to} \lambda$ by Theorem \ref{T:tsukada}(2).

It remains to construct the conditional expectation.  Since $s(P) = J(1)$ by assumption, we only need to consider elements in $J(1)\mathcal{M}_2 J(1)$.  We use $J(1) = \pi(1) + \pi'(1)$ to write a generic element as
$$y = \left( \begin{smallmatrix} a & b \\ c & d \end{smallmatrix} \right),$$
where $a \in \pi(1) \mathcal{M}_2 \pi(1)$, etc.  By Theorem \ref{T:tsukada}(2) we have
\begin{equation} \label{E:slim}
P(y) = s-\lim \bar{E}_\alpha \circ P(y) = s-\lim S_{\lambda_\alpha} F_\alpha(y).
\end{equation}
Since $F_\alpha$ is a conditional expectation and $J(\mathcal{M}_\alpha)$ is inside the diagonal $J(\mathcal{M}_1)''$, we may write
$$F_\alpha \left( \begin{smallmatrix} a & b \\ c & d \end{smallmatrix} \right) = \left( \begin{smallmatrix} F^1_\alpha(a) & 0 \\ 0 & F^2_\alpha(d) \end{smallmatrix} \right).$$
Then \eqref{E:slim} becomes
\begin{equation} \label{E:longp}
P(y) =
\end{equation}
$$ s-\lim \left( \begin{smallmatrix} \pi(\lambda_\alpha) F^1_\alpha(a) + \pi(1-\lambda_\alpha) [\pi \circ \pi'^{-1}(F^2_\alpha(d))] & 0 \\ 0 & \pi'(\lambda_\alpha) [\pi' \circ \pi^{-1}(F^1_\alpha(a))] + \pi'(1-\lambda_\alpha) F^2_\alpha(d)  \end{smallmatrix} \right),$$
which we claim can be written as
\begin{equation} \label{E:pfactor}
P(y) = S_\lambda \circ F(y), \text{ where } F(y) = s-\lim \left( \begin{smallmatrix} F^1_\alpha(a) & 0 \\ 0 & F^2_\alpha(d)  \end{smallmatrix} \right).
\end{equation}
To establish the claim, we check that (i) the strong limits in the definition of $F$ do exist, (ii) the factorization \eqref{E:pfactor} is correct, (iii) $F$ is normal, (iv) $F$ is contractive, (v) the range of $F$ is contained in $J(\mathcal{M}_1)''$, and (vi) $F$ fixes $J(\mathcal{M}_1)''$.

We have from the above that
$$P \left( \left( \begin{smallmatrix} a & 0 \\ 0 & 0 \end{smallmatrix} \right)
\right) = s-\lim \left( \begin{smallmatrix} \pi(\lambda_\alpha) F^1_\alpha(a) & 0 \\ 0 & \pi'(\lambda_\alpha) [\pi' \circ \pi^{-1}(F^1_\alpha(a))] \end{smallmatrix} \right),$$
so $\pi(\lambda_\alpha) F^1_\alpha(a)$ is strongly convergent.   We now represent $\pi(\mathcal{M})$ faithfully and nondegenerately as an algebra of operators on a Hilbert space, and let $e_n = e[1/n, 1]$ be spectral projections of $\pi(\lambda)$.  Since $\lambda$ is nonsingular, these projections increase to $\pi(1)$.  Given any vector $\xi$, we have
\begin{align*}
\| [F^1_\alpha(a) - F^1_\beta(a)]  \xi \| &\sim \| [ F^1_\alpha(a) - F^1_\beta(a)] e_n \xi \| \\
&= \| [ F^1_\alpha(a) - F^1_\beta(a)] \pi(\lambda) (\pi(\lambda)^{-1} e_n \xi) \| \\
&= \| (( F^1_\alpha(a)[\pi(\lambda) - \pi(\lambda_\alpha)] + [\pi(\lambda_\alpha) F^1_\alpha(a) - \pi(\lambda_\beta) F^1_\beta(a)] \\
&\qquad + F^1_\beta(a)[\pi(\lambda) - \pi(\lambda_\beta)] ) (\pi(\lambda)^{-1} e_n \xi) \| \\
&\to 0 \text{ as $\alpha, \beta$ increase.}
\end{align*}
Note that the first approximation can be made independent of $\alpha$ and $\beta$ since $\{F^1_\alpha(a)\}$ is a norm-bounded set.  The same argument establishes the strong convergence of $F^2_\alpha(d)$, and (i) is obtained.  Since $\pi(\lambda_\alpha), F^1_\alpha(a), F^2_\alpha(d)$ are strongly convergent, (ii) follows from \eqref{E:longp}.

For (iii), again it suffices by symmetry to check that 
\begin{equation} \label{E:fnormal}
\pi(1) \mathcal{M}_2 \pi(1) \ni x_\gamma \overset{s}{\to} x \Rightarrow F \left( \left( \begin{smallmatrix} x_\gamma & 0 \\ 0 & 0 \end{smallmatrix} \right) \right) \overset{s}{\to} F \left( \left( \begin{smallmatrix} x & 0 \\ 0 & 0 \end{smallmatrix} \right) \right).
\end{equation}
(All limits in this paragraph are along increasing $\gamma$.)  Since $P$ is normal, we have
$$ P \left( \left( \begin{smallmatrix} x_\gamma & 0 \\ 0 & 0 \end{smallmatrix} \right) \right) \overset{s}{\to} P \left( \left( \begin{smallmatrix} x & 0 \\ 0 & 0 \end{smallmatrix} \right) \right)$$
and then
$$ S_\lambda \left( \left( \begin{smallmatrix} s-\lim_\alpha F^1_\alpha (x_\gamma) & 0 \\ 0 & 0 \end{smallmatrix} \right) \right) \overset{s}{\to} S_\lambda \left( \left( \begin{smallmatrix} s-\lim_\alpha F^1_\alpha(x) & 0 \\ 0 & 0 \end{smallmatrix} \right) \right).$$
Reading off the upper left entry,
$$\pi(\lambda) [s-\lim_\alpha F^1_\alpha (x_\gamma)] \overset{s}{\to} \pi(\lambda) [s-\lim_\alpha F^1_\alpha (x)].$$
By an approximation argument similar to that of the previous paragraph,
$$s-\lim_\alpha F^1_\alpha (x_\gamma) \overset{s}{\to} s-\lim_\alpha F^1_\alpha (x).$$
But this is exactly the conclusion of \eqref{E:fnormal}, and we have (iii).

Now (iv) and (v) are automatic from the form of $F$, and (vi) follows from the normality of $F$ and the fact that $F$ fixes $\cup J(\mathcal{M}_\alpha)''$.  The proof is complete.

\end{proof}

The preceding propositions may seem somewhat technical, but they can be summarized nicely.

\begin{theorem} \label{T:astotyp}
Let $J:\mathcal{M}_1 \to \mathcal{M}_2$ be a normal Jordan *-monomorphism, and let $P: \mathcal{M}_2 \to J(\mathcal{M}_1)$ be a normal positive projection.  Each condition below implies its successor:
\begin{enumerate}
\item $\mathcal{M}_1$ is AS;
\item the projection $P$ factors through a conditional expectation onto $J(\mathcal{M}_1)''$;
\item the map $\varphi^{1/p} \mapsto (\varphi \circ J^{-1} \circ P)^{1/p}$, $\varphi \in (\mathcal{M}_1)_*^+$, extends linearly to an isometry from $L^p(\mathcal{M}_1)$ to $L^p(\mathcal{M}_2)$.\
\end{enumerate}
\end{theorem} 

\section{Conclusion}

From Theorems \ref{T:main} and \ref{T:astotyp}, we have arrived at a complete description of all the isometries from $L^p(\mathcal{M}_1)$ to $L^p(\mathcal{M}_2)$, provided that $\mathcal{M}_1$ is AS.  They are typical, arising via $\eqref{E:typical}$, and arbitrary data $J, P, w$ are allowed.  When $\mathcal{M}_1$ is only assumed to have EP$p$ and EP1, which may be weaker, we can still say that all $L^p$ isometries are typical.

The main motivation for this paper was to develop Watanabe's ideas, especially EP, as far as possible.  This certainly provides new methods and information, but we have only been able to apply them to Questions \ref{T:conj} and \ref{T:conj2} when the initial algebra is well-approximated by semifinite algebras.  We do not claim that further work in this direction is \textit{necessary} for a solution to one or both of these questions.  It certainly seems possible that a new technique may produce a relatively straightforward solution - for example, the paper [S1] handles the surjective case by different methods.  Nonetheless, we do think that EP$p$, AS, and condition (2) from Theorem \ref{T:astotyp} are worth investigating on their own merits, and our work here shows their relation to $L^p$ isometries.

\end{document}